\documentclass[a4paper,11pt]{article}

\usepackage[pdftex]{hyperref}
\usepackage{stmaryrd}
\usepackage[normalem]{ulem}
\usepackage{graphicx}
\usepackage{amsfonts}
\usepackage{amsmath,amssymb,mathrsfs,amsthm}
\usepackage[all]{xy}
\usepackage{color}
\usepackage{times}
\usepackage{enumerate,vmargin}
\usepackage{verbatim}
\usepackage{paralist} 
\usepackage[labelfont=bf,font=it]{caption}
\usepackage[numbers,sort&compress]{natbib}

\setpapersize{A4}
\setmarginsrb{2.5cm}{2cm}{2.5cm}{2cm}{.4cm}{4mm}{0cm}{7.5mm}

\hypersetup{
    unicode      = false,     
    pdftoolbar   = true,      
    pdfmenubar   = true,      
    pdffitwindow = true,      
    pdfnewwindow = true,      
    colorlinks   = true,      
    linkcolor    = blue,      
    citecolor    = blue,      
    filecolor    = blue,      
    urlcolor     = blue       
}

\theoremstyle{plain}
\newtheorem{thm}{Theorem}

\theoremstyle{plain}
\newtheorem{lem}[thm]{Lemma}
\newtheorem{prop}[thm]{Proposition}

\theoremstyle{definition}

\newtheorem{rem}{Remark}

\newcommand{\noi}{\noindent}

\definecolor{grey}{rgb}{.5,.5,.5}
\definecolor{cnic}{rgb}{.9,0.5,0}
\definecolor{mycolor}{rgb}{0,.6,0.4}

\newcommand{\cT}{\mathcal T}
\newcommand{\eqd}{\stackrel{\mathrm{(d)}}=}

\newcommand{\bbT}{\mathbb T}

\newcommand{\rT}{\mathrm T}
\newcommand{\N}{\mathbb N}

\newcommand{\bN}{\mathbf N}

\newcommand{\cP}{\mathcal P}

\newcommand{\bQ}{\mathbf Q}

\newcommand{\llb}{\llbracket}
\newcommand{\rrb}{\rrbracket}
\newcommand{\degree}{\mathtt{deg}}
\newcommand{\diist}{\mathtt{dist}}
\newcommand{\diam}{D}

\newcommand{\bbU}{\mathbb{U}}
\newcommand{\cF}{\mathcal{F}}
\newcommand{\bP}{\mathbf{P}}
\def\cq{$\hfill \square$}
\def\cqfd{$\hfill \blacksquare$}
\newcommand{\bC}{\mathbf{C}}
\newcommand{\bbR}{\mathbb{R}}

\newcommand{\bE}{\mathbf{E}}
\newcommand{\lam}{\lambda}
\newcommand{\ep}{\varepsilon}
\newcommand{\un}{\mathbf{1}}
\newcommand{\indi}{\un}
\newcommand{\cI}{\mathcal I}
\newcommand{\bep}{\mathcal{E}}

\newcommand{\Trp}{\mathrm{T}^{_{(p)}}_{^{\lfloor rp\rfloor}}}
\newcommand{\epp}{\mathcal{E}^{_{(p)}}_{^{\lfloor rp\rfloor}}}

\newcommand{\Tpl}{\bbT_{\! {\rm pl}}}
\newcommand{\Tor}{\bbT_{\! {\rm or}}}
\newcommand{\Wt}{\mathrm{Wt}_{ \mu}}

\title{ \textsc{Scaling limits for a family of unrooted trees}}
\author{Minmin \textsc{Wang}
\thanks{Conicet-UBA. Departamento de Matemática, Universidad de Buenos Aires,
C1428EGA, Buenos Aires - Argentina. Email: wangminmin03@gmail.com}
}

\date{July 2016}

\begin{document} 

\maketitle

\begin{abstract}
We introduce weights on the unrooted unlabelled plane trees as follows: let $\mu$ be a probability measure on the set of nonnegative integers whose mean is bounded by $1$; then the $\mu$-weight of a 
plane tree $t$ is defined as $\Pi \, \mu (\mathtt{degree} (v) \! -\!1)$, where the product is over the set of vertices $v$ of $t$. We study the random plane tree with a fixed diameter $p$ sampled according to probabilities  
proportional to these $\mu$-weights and we prove that, under the assumption that the sequence of laws 
$\mu_p$, $p\! \geq \! 1$, belongs to the domain of attraction of an infinitely divisible law, the scaling limits of such random plane trees are random compact real trees called 
the unrooted L\'evy trees,  which  have been introduced in \cite{DuWa14}.

\smallskip


\smallskip

\end{abstract}

\section{Introduction}
\label{introsec}
Models of random \textit{rooted} trees have been extensively studied (see for instance Aldous \cite{aldcrt2}, Devroye \cite{Devroye}, Le Gall \cite{LGCornell}, Drmota \cite{Drmotabook}, Janson \cite{Jansonsurvey}) 
often because of their connections with branching processes: an eminent example being the model of Galton--Watson trees. 
The scaling limits of Galton--Watson trees are L\'evy trees, introduced by Le Gall and Le Jan \cite{legalllejan}. L\'evy trees extend Aldous' notion of Brownian Continuum Random Tree \cite{aldcrt1, aldcrt3}; they 
describe the genealogy of continuous-state branching processes; they are closely related to fragmentation and coalescent processes: see Miermont \cite{Mi03, Miermont05}, Haas \& Miermont \cite{HaMi04}, Goldschmidt \& Haas \cite{GolHaa}, Abraham \& Delmas \cite{AbDel08, AbDe15}. 
Their probabilistic and fractal properties are studied in Duquesne \& Le Gall \cite{Duquesne02, Duquesne05}. 

However, there are situations where unrooted trees arise naturally. In this work, we focus on \emph{unrooted and unlabelled plane trees} that are graph-trees embedded into the oriented plane and considered up to orientation-preserving homeomorphisms (see Section \ref{sec: notation} for a precise definition). 

We propose here a model of random plane trees defined as follows. 
Let $\mu$ be a probability measure on $\mathbb N:=\{0, 1, 2, \dots\}$. 
The \textit{$\mu$-weight} of a plane tree $t$ is the quantity: 
\begin{equation}\label{def: Wt}
W_\mu(t):=\prod_{v \text{ vertex of } t} \mu (\degree (v) \! -\!1),
\end{equation}
where $\degree(v)$ denotes the degree of the vertex $v$ in $t$. These $\mu$-weights induce a  probability measure on the sets of plane trees with a fixed diameter. More precisely, the \emph{graph distance} between two vertices of a tree $t$ is the number of edges on the unique path of $t$ between the two  vertices. Then the \emph{diameter} of $t$ is the maximum distance between any pair of vertices of $t$. For $p\ge 0$, let $\mathrm T_p$ be a random plane tree with diameter $p$ 
such that the probability of the event $\mathrm T_p=t$ is proportional to $W_\mu(t)$, for all plane tree $t$ with diameter $p$ (see Section \ref{sec: result} below for a more careful definition of $\mathrm T_p$). 

The above definition of the random plane tree $\mathrm T_p$ is inspired by simply generated trees introduced by Meir \& Moon (\cite{MM_simply}). 
Our study of $\mathrm T_p$ is, on the other hand, motivated by a previous work \cite{DuWa14}, where distributional properties of the diameters of general L\'evy trees are examined. In particular, from the study there, a notion of \emph{unrooted L\'evy trees} naturally arises. 
Intuitively, an unrooted  L\'evy tree with diameter $r$ is obtained from two independent L\'evy trees conditioned to have height $r/2$ by connecting their roots. It is shown in \cite{DuWa14} that 
the spinal decomposition of an unrooted L\'evy tree along its longest geodesic exhibits a remarkably simple form. On the other hand, a classical (rooted) L\'evy tree can be obtained from the unrooted one by picking a uniform point as the root. (See \cite{DuWa14} or Section \ref{sec: Levy} for the precise statements.)  
Because of all these properties, the model of unrooted L\'evy trees seems to us an interesting object to study. 

The work \cite{DuWa14} deals with the continuum trees. Here, we consider the discrete counterpart (namely, $\mathrm T_p$) and  the convergence of the discrete trees when the diameters tend to infinity. In the main result (Theorem \ref{main}, see also Remark \ref{rem: main} there), we show that unrooted L\'evy trees appear in the scaling limits of $\mathrm T_p$. 
This result could be viewed as the analog for unrooted plane trees to the Duquesne--Le Gall's Theorem \cite{Duquesne02}, which establishes the convergence of rescaled Galton--Watson trees to L\'evy trees. 

As an essential ingredient of the main proof, we also prove a limit theorem for Galton-Watson trees conditioned to have a fixed  large height (Proposition \ref{bdkjfvs}), which might be of independent interest. This result extends a previous one due to Le Gall \cite{LeG10} in the Brownian case, which has been proved by a different method. 
The idea here is to perform a simple transformation on Galton--Watson trees, which consists in extracting a Galton--Watson tree conditioned to have a height $k$ from a Galton--Watson tree conditioned to have a height $\ge k$. 

There are many possible ways to condition a random tree to be large. Previous works (see for example Aldous \cite{aldcrt3}, Marckert \& Miermont \cite{MaMi11}, Haas \& Miermont \cite{HaMi12}) mainly focus on random trees with a fixed progeny. Here, we condition trees by their heights or by their diameters, which are more adapted to the limit objects considered here.  We refer to the paper \cite{LeG10} of Le Gall for a discussion on various conditionings and their connections with the excursion measures. 

Another feather of the current work is that, unlike most of the previous works on the limit theorems of random trees, we look at unrooted trees rather than rooted ones. Due to their internal symmetries, unrooted trees turn out to be more difficult to deal with. Here, we employ the centers of a tree: 
our proof below of Theorem \ref{main} relies on a decomposition of $\mathrm T_p$ at its centers. 
After handling a technical point on the central symmetries, we show that asymptotically, this decomposition results in  two independent rooted trees with fixed heights. 

\medskip
The paper is organized as follows. In Section \ref{sec: notation}, we introduce the necessary notation  and provide background on different classes of trees that are studied here. The main results, Proposition \ref{bdkjfvs} and Theorem \ref{main},  are stated in Section \ref{sec: result}. Their proofs are found in Sections \ref{Pfbdkjfvs} and \ref{mainthpf}. 

\section{Preliminaries and notation}
\label{sec: notation}

Before stating the main results of the paper, we recall here some notation and 
the definitions of the various classes of discrete trees that appear in the proofs of the main results. We wish to emphasize on differences and connections between plane trees and ordered rooted trees. We will see that equivalent classes of edge-rooted plane trees correspond to ordered rooted trees. 
On the other hand, we need to take into account the symmetry of a plane tree when performing a rooting of the tree. 
We also give a brief introduction to real trees and L\'evy trees, upon which the construction of the limit objects relies.  
For a more extensive account on discrete trees, we refer to Drmota's book \cite{Drmotabook}; for the technical background on L\'evy trees we refer to Duquesne \& Le Gall \cite{Duquesne02, Duquesne05}; see also Evans \cite{evans05} and the references there for more information on real trees. 

Unless otherwise specified, all the random variables that we mention here are defined on a common probability space $(\Omega, \cF, \bP)$. 

\subsection{Diameter and centers of trees} 
\label{sec: ctr}
A \emph{tree} is a connected graph without cycle. We only consider \textit{finite} trees. 
Two vertices $v, w$ are \emph{adjacent} if there is an edge between them. In that case, we write $v\sim w$; we also use the notation $(v, w)$ to indicate the \emph{oriented} edge pointing from $v$ to $w$. 
The \emph{degree} of a vertex $v$ in a tree $t$ is the number of its adjacent vertices: $\degree (v)= |\{ \textrm{$w$ vertex: $v\sim w$} \} |$. 
The \emph{size} of a tree $t$, denoted by $|t|$, is the number of its vertices. 
A \textit{path} from $v$ to $w$ in the tree $t$ is a sequence of adjacent vertices $v=v_0\sim v_1 \sim \cdots \sim v_n=w $. We denote by $\llb v, w\rrb $ the unique self-avoiding path joining $v$ to $w$. Then 
the \emph{graph distance} between $v$ and $w$, denoted by $\diist(v,w)$, is the number of edges of $\llb v, w\rrb$, which we also refer to as the \emph{length} of the path.

For a tree $t$, we denote respectively by $V$ and $E$ the sets of its vertices and of its edges. 
For all $v\in V$, we set 
\[
\Gamma (v,t)= \max_{w\in V} \; \diist (v, w).
\]
We then define respectively
\[
\diam(t)=\max_{v\in V} \; \Gamma (v,t) \quad \textrm{and} \quad \gamma (t) = \min_{v\in V} \; \Gamma (v,t)\; .
\]
We say that $ \diam(t)=\max_{v,w\in V} \diist(v,w)$ is the \emph{diameter} of the tree $t$. The following notion plays an important role in this work. 
\begin{equation}
\label{centerdef}
 \textrm{A vertex $v\in V$ is a \textit{center} of the tree $t$ if $ \Gamma (v,t)= \gamma (t)$.}
 \end{equation}
In  \cite{Jor1869} (1869) C.~Jordan  proved that a tree has either one or two centers. More precisely, for the tree $t$, we have a dichotomy in the number of centers of $t$ depending on the parity of its diameter. 

\smallskip
\begin{compactenum}
\item[--] \textit{Bi-centered case:} if $\diam(t)$ is odd, then $\gamma (t)= \frac{1}{2} (\diam (t)+1)$ and there are two adjacent centers $c, c^\prime$. Moreover in this case, for any path $\llb v, w\rrb$ such that 
$\diist(v,w)= \diam(t)$, $c$ and  $c^\prime$ are the two midpoints of the path. Namely, they are the
only vertices in $\llb v, w\rrb$ at distance $\gamma(t)$ from either $v$ or $w$.
\item[--] \textit{Uni-centered case:} if $\diam(t)$ is even, then $\gamma (t)= \frac{1}{2} \diam (t)$ and there is a unique center $c$. Moreover in this case, for any path $\llb v, w\rrb$ such that 
$\diist(v,w)= \diam(t)$, $c$ belongs to the path and  is situated at equal distance from $v$ and $w$. 
\end{compactenum} 
As it turns out, in our study of unrooted  trees, centers are convenient choices for roots.

\subsection{Ordered rooted trees and Galton-Watson trees}
\label{sec: or}

Let us first recall Ulam's coding for ordered rooted trees. To this end, write $\N^*\! = \! \{ 1, 2, \ldots \}$. Define $\bbU\!:= \{ \varnothing \}\cup  \bigcup_{n\geq 1} (\N^*)^n$ so that an element $u$ of $\bbU$ is a finite sequence: $u=(a_1, \dots, a_n)$ for some integer $n\ge 0$ and $a_i\in \N^*$, $1\le i\le n$. We say that a finite subset $s$ of $\bbU$ is an \emph{ordered rooted tree} if it satisfies the following three properties.
\begin{compactenum}
\item[(a)]  $\varnothing \in s$.
\item[(b)] If $u= (a_1, \ldots, a_n) \in s \backslash \{ \varnothing \}$, then $(a_1, \ldots, a_{n-1}) \in s$. 
\item[(c)] For all $u= (a_1, \ldots, a_n) \in s $, there exists a nonnegative integer $k_u(s)$ such that if $k_u(s) \geq 1$, then $(a_1, \ldots, a_n, a) \in t$, for all $1\leq a \leq k_u(s)$. 
\end{compactenum}
Alternatively, an ordered rooted tree $s$ can be identified as a family tree: $\varnothing$ is the common ancestor; the elements $(1), (2), \dots, (k_\varnothing(s))$ form the first generation, ranked in their birth orders; more generally, an individual $u=(a_1, a_2, \cdots, a_n)\in s$ has exactly $k_u(s)$ children, namely $(a_1, a_2, \cdots, a_n, j)$, $j=1, 2, \dots, k_u(s)$. Observe that the ancestors of $u$ are given by  $(a_1, a_2, \dots, a_i)$, $0\le i\le n$.  

We denote by $\Tor$ the set of ordered rooted trees. 
Let us recall a well-known result on the enumeration of ordered rooted trees. For each $n\in \mathbb N$, the number of ordered rooted trees of size $n\!+\!1$ is given by $C_n=\frac{1}{n+1} \binom{2n}{n}$, the $n$-th Catalan number. 

\bigskip

Let $\mu= (\mu(k))_{k\geq 0}$ be a probability distribution on  
$\N\! :=\!  \{ 0, 1, \ldots \}$. Since we only consider finite trees, let us assume that the mean of $\mu$ is bounded by $1$. We also suppose that the support of $\mu $ is not contained in the set $\{0, 1\}$ to exclude trivial cases.  
We summarize our assumptions  on $\mu$ as follows: 
\begin{equation}
\label{subcrGW}
 \forall k\in \N, \quad \mu(k) \geq 0, \quad \sum_{k\in \N} \mu (k)= 1 \, ,  \quad \sum_{k\in \N} k \mu (k)  \leq 1 \quad \textrm{and} \quad \exists\, k \geq 2 : \mu (k) >0 \; .
\end{equation} 
For each finite ordered rooted tree $s\in \Tor$, we set 
\begin{equation}
\label{GWlawdef}
 P_\mu (s)=  \prod_{u\in s}  \mu\big(k_u(s)\big) .
\end{equation} 
Standard arguments (see for example Neveu \cite{Neveu}) show that $P_\mu$ defines a probability measure on $\bbT_{\! {\rm or}}$; $P_\mu$ is then called  \textit{Galton--Watson law with offspring distribution $\mu$}. In what follows, a GW($\mu$)-tree refers to a random variable, say $\tau$, defined on the probability space $(\Omega, \cF, \bP)$ taking values in $\Tor$ such that $\bP(\tau=s)\! = \! P_\mu (s)$, for all $s\in \Tor$. 

Galton--Watson trees are known to be closely related to 
simply generated trees; we refer to Aldous \cite{aldcrt2}, Devroye \cite{Devroye}, Drmota \cite{Drmotabook}, Janson \cite{Jansonsurvey} and Kennedy \cite{Kennedy_cond} for a more extensive account on this subject.

\subsection{Plane trees } 
\label{sec: pl}
Plane trees are graph-trees embedded into the oriented plane $\cP$, considered up to orientation preserving homeomorphisms from $\cP$ to $\cP$. 
More precisely, an \textit{oriented edge} in $\cP$ (distinct from a loop) is a continuous and injective function $\varepsilon \!  : \! [0, 1] \! \rightarrow \! \cP$, considered up to reparametrization (i.e.\ strictly increasing homeomorphisms from $[0, 1]$ to $[0, 1]$). 
The \textit{tail} of the oriented edge $\varepsilon$ is $\varepsilon(0)$ and its \textit{head} is $\varepsilon(1)$. The \textit{reversed} edge is $\overline{\varepsilon} (t)=\varepsilon(1-t)$, $t\in [0, 1]$, also considered up to reparametrization. A non-oriented \textit{edge} is then given by $e\! =\! \{ \varepsilon, \overline{\varepsilon}\}$. The \textit{endpoints} of $e$ are $\{ \varepsilon (0), \varepsilon(1)\}$ and 
its \textit{inner part} is $\varepsilon ((0, 1))$. Note that the endpoints and the inner part of an edge do not depend on any particular parametrization.  
An \emph{embedded tree in the plane} is a pair $t=(V, E)$ such that
\vspace{1mm}
\begin{compactenum}
\item[(a)] $V\subset \cP$ is finite;
\item[(b)] $E$ is a connected subset of $\cP$ formed by a finite set of edges (as defined above) whose inner parts are pairwise disjoint and do not intersect $V$ and whose endpoints belong to $V$; 
\item[(c)] $|V| \! =\! |E|+1$.       
\end{compactenum}
\vspace{0.5mm}
Two embedded plane trees $t\! =\! (V, E)$ and $t^\prime \! =\! (V^\prime , E^\prime)$ are said to be \emph{equivalent} if there exists an orientation preserving homeomorphism $h \! :\!  \cP \! \rightarrow \! \cP$ such that $V^\prime\! = \! h(V)$ and $E^\prime \! =\!  \{\{ h\! \circ \! \varepsilon, h \! \circ\!  \overline{\varepsilon} \}: \{\varepsilon, \overline{\varepsilon}\}\in E\}$. The equivalence class of an embedded plane tree is referred to as a \emph{plane tree} in the rest of the paper. We denote by $\bbT_{\! {\rm pl}}$ the set of plane trees. 

\begin{figure}[tp]
\centering
\includegraphics[scale=.7]{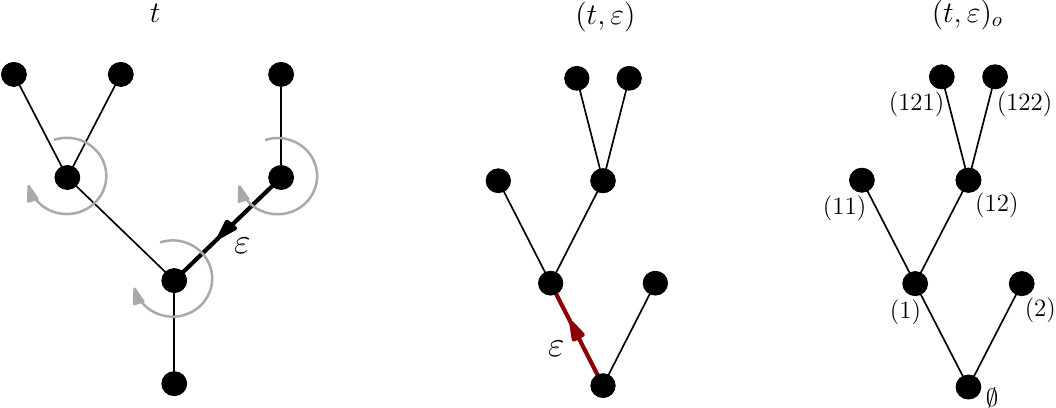}
\caption{\label{fig: contour} On the left, an embedded plane tree $t$: the plane orientation induces a  cyclic order on the neighbors of each vertex.  In the middle, the edge-rooted plane tree $(t, \varepsilon)$ obtained by rooting $t$ at the oriented edge $\ep$.
On the right, the ordered rooted tree $(t, \ep)_o$ associated with $(t, \ep)$. }
\end{figure}

A plane tree $t=(V, E)$ naturally induces a graph-tree; but it also carries additional structures inherited from the oriented plane. In particular, for each vertice $v \in V$, the orientation of the plane induces a cyclic order on the set of vertices that are adjacent to $v$ (namely, the set of the \textit{neighbors of $v$}); see Figure \ref{fig: contour}. This notion of cyclic order will be useful in what follows.

In \cite{HarPriTut64}, Harary, Prins \& Tutte deduce a 
functional equation for the generating function of the numbers of plane trees with given sizes, which eventually leads to the following closed formula due to Walkup \cite{Wal72}: for each $n\in \mathbb N$, 
\begin{equation}
\label{enumpla}
\big| \{ t \! \in \! \bbT_{\! {\rm pl}}\!: |t| \! =\! n+1 \} \big| \! =\!   \frac{1}{2n(n+1)} \binom{2n}{n} + \frac{1}{4n} \binom{n+1}{\frac{n+1}{2}}\un_{\{n \text{ is odd}\}}+ \frac{\varphi (n)}{n}+ \frac{1}{2n} \!\!\! \sum_{\substack{d | n \\ 1<d<n}} \!\!\! \varphi \Big(\frac{n}{d}\Big) \binom{2d}{d}, 
\end{equation}
where $\varphi$ stands for Euler's totient function. The somewhat complicated form of \eqref{enumpla} is an indication of  the presence of internal symmetry in plane trees; see Remark \ref{internalsym} below. 

Let us mention that plane trees are particular instances of planar maps: they are planar maps with one face. We refer to Mohar \& Thomassen's book \cite{MohTho01} for a more precise account on embedded graphs on surfaces and to Lando \& Zvonkin's book \cite{LanZvo04} for a combinatorial definition of planar maps.

\paragraph{Edge-rooted plane trees and ordered rooted trees. } 
A plane tree with a distinguished oriented edge is called an \textit{edge-rooted plane tree}. Two edge-rooted plane trees $(t, \varepsilon)$ and $(t^\prime, \varepsilon^\prime)$
 are said to be \emph{equivalent} if there exists an orientation preserving homeomorphism $h: t\to t^\prime$ satisfying $\varepsilon^\prime \! =\!  h \! \circ \! \varepsilon$. 
Let $(t, \varepsilon)$ be an edge-rooted plane tree. We can associate with it an ordered rooted tree in the following way (see also Fig.~\ref{fig: contour}). Let us employ the terminology of family tree 
and recall that there is a cyclic order on the neighbors of each vertex of $t$ which is induced by the orientation of the plane. Let $\rho\!:=\varepsilon(0)$, the tail of $\varepsilon$. We view $\rho$ as the common ancestor. Let $v_1, v_2, \dots, v_{\deg(\rho)}$ be the neighbors of $\rho$ ordered in such a way that $v_1=\varepsilon(1)$ and that $v_i$ is next to $v_{i-1}$ in the cyclic order, for all $2\le i\le \deg(\rho)$. Then the sequence $(v_i)_{1\le i\le \deg(\rho)}$ forms the first generation ranked in their birth orders. 
More generally, for a vertex $u\ne \rho$, let us write $v_0, v_1, \dots, v_{\deg(u)-1}$ for its neighbors ordered in such a way that $v_0$ is the unique vertex of $\llb\varepsilon(0), u\rrb$ that is adjacent to $u$ and that $v_i$ is next to $v_{i-1}$ in the cyclic order, for all $1\le i\le \deg(u)-1$. Then $v_1, \dots, v_{\deg(u)-1}$ are the children of $u$, ranked in their birth orders, while $v_0$ is the parent of $u$. 
In this way, we can readily associate with $(t, \ep)$ an ordered rooted tree which we denote by $(t, \ep)_o\in \Tor$. It is straightforward to check that if $(t, \ep)$ and $(t', \ep')$ are two equivalent edge-rooted plane trees, then $(t, \ep)_o=(t', \ep')_o$. On the other hand, for an ordered rooted tree $s\in \Tor$ satisfying $|s|>1$, we can always find an embedded plane tree $t$ and an oriented edge $\ep$ of $t$ such that $(t, \ep)_o \! =\!  s$. To sum up, \textit{there is a bijection between the set of equivalence classes of edge-rooted plane trees and the set of ordered rooted trees with size $>1$}.

\begin{figure}[htp]
\centering
\includegraphics[scale=.6]{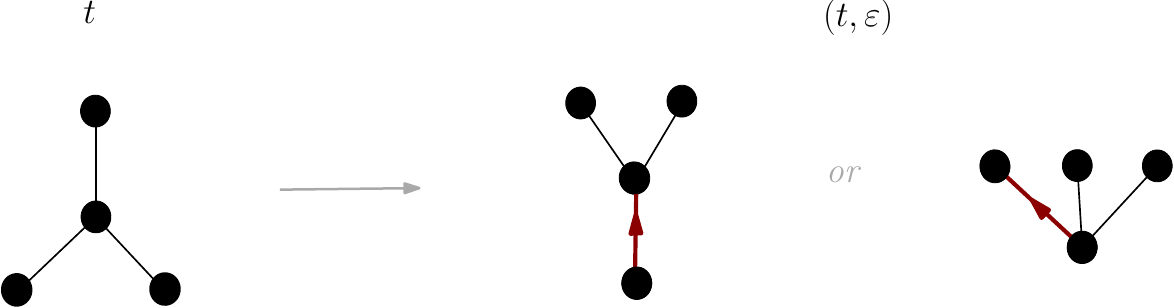}
\caption{\label{fig: symmetry} On the left, a plane tree $t$ which has $6$ oriented edges. By rooting it on each of these edges, we obtain $2$ equivalent classes of edge-rooted trees,  as illustrated on the right.  }
\end{figure}

At this point, let us make an important remark on the number of possible ways in rooting a plane tree, which turns out to be one of the technical points that we need to deal with in the proof of the main theorem. 
\begin{rem}
\label{internalsym} 
Due to a potential internal symmetry of the tree, the mapping $(t, \ep) \mapsto (t,\ep)_o$ is surjective but not injective in general, since different choices of edges may give rise to equivalent edge-rooted plane trees; 
see Figure \ref{fig: symmetry} for an example.  Indeed, let us observe from \eqref{enumpla} that for all $n\! \geq \! 1$, we have 
\[
2n \cdot \big| \big\{ t \in \Tpl: |t| =  n+1 \big\} \big|  >  \big| \big\{ t  \in  \Tor : |t|  = n+1 \big\} \big|=\frac{1}{n+1} \binom{2n}{n} \,.
\]
On the other hand, note that 
\[
2n \cdot \big| \big\{ t  \in  \Tpl: |t|  = n+1 \big\} \big| \ \overset{n\rightarrow \infty}{\sim} \ \big| \big\{ t  \in  \Tor : |t|  =n+1 \big\} \big| \; , 
\]
which suggests that a ``typical'' large plane tree has no internal symmetries. \cq 
\end{rem}

\begin{figure}[tp]
\centering
\includegraphics[scale=.7]{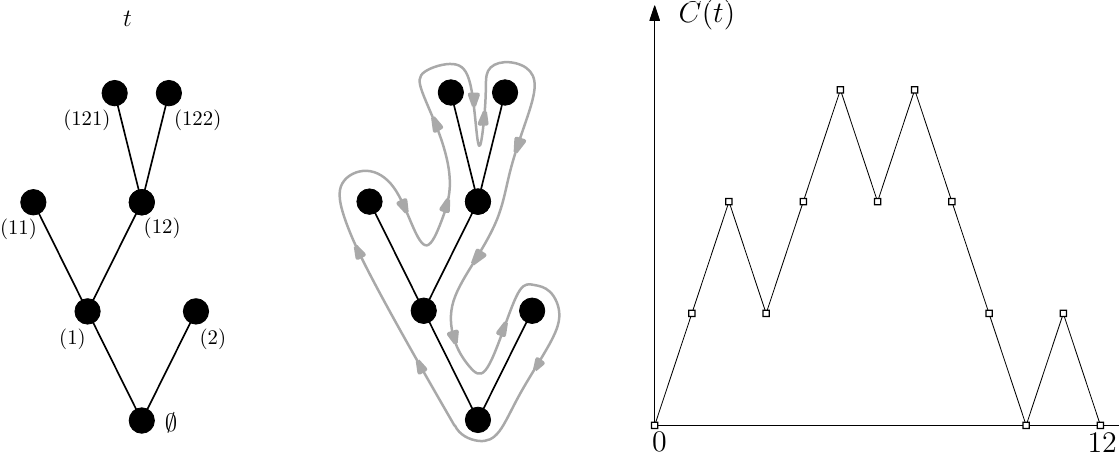}
\caption{\label{fig: cont} The contour function of an ordered rooted tree $t$ (on the left) with $6$ edges. 
In the middle, the grey line illustrates the exploration of the particle.
On the right, its contour function $C(t)$. }
\end{figure}

\paragraph{Contour functions of edge-rooted plane trees and ordered rooted trees. } 
We will use contour functions to study trees, whose definition we recall here. 
In the first place, let $t$ be an ordered rooted tree embedded into the plane in such a way that the common ancestor 
$\varnothing$ is located at the origin and 
the children of each vertex appear from left to right in increasing order of their birth orders. In a less formal way, imagine that each edge is a line segment of length $1$ and a particle explores the embedded tree with speed $1$ from the common ancestor, in a left-to-right way, backtracking as least as possible. It terminates its exploration at time $2(|t|-1)$, during which each edge is visited exactly twice by the particle. 
Denote by $C_s(t)$ the distance of the particle from $\varnothing$ at time $s$.  Then the continuous function $(C_s(t), 0\le s\le 2(|t|-1))$ is called the \emph{contour function of the ordered rooted tree $t$}. 
See Fig.~\ref{fig: cont}. We refer to Duquesne \cite{Duquesne03} for a formal definition. 
Note that $t$ is characterized by its contour function. In particular, the graph distance $\diist$ of $t$ can be found as follows. For $p=0, 1, \dots, 2(|t|-1)$, denote by $v_p$ the vertex of $t$ visited by the particle at time $p$. 
Then, for all integers $p, q$ such that $0\leq p \leq q \leq 2(|t|-1)$, we have 
\[
\diist(v_p, v_q) = C_p(t) + C_q(t) -  2 \inf_{s\in [p, q]}  C_s(t) \; .
\]
Next, we recall from above that an edge-rooted plane tree $(t, \ep)$ can be associated with an ordered rooted tree $(t, \ep)_o$. Then we define the \emph{contour function of the edge-rooted tree $(t, \ep)$} to be that of $(t, \ep)_o$ and we denote it by $C(t, \ep)$.

\subsection{L\'evy trees} 
\label{sec: Levy}
The class of \textit{L\'evy trees} is introduced by Le Gall \& Le Jan \cite{legalllejan} (see also Duquesne \& Le Gall \cite{Duquesne02}). It 
extends Aldous' \textit{Continuum Random Tree} \cite{aldcrt1,aldcrt2,aldcrt3}, which will sometimes be referred to as the Brownian case. L\'evy trees are random compact metric spaces, or more specifically, random compact real trees. Let us first recall the definition of real trees.

\paragraph{Real trees. } One way to generalize graph-trees is to consider geodesic metric spaces without loops. More precisely, a metric space $(T, \delta)$ is a \emph{real tree} if the following two properties hold.
\vspace{1mm}
\begin{compactenum}[(i)]
\item For all $\sigma, \sigma^\prime \in T$, there is a unique isometric mapping $q: [0, \delta(\sigma, \sigma^\prime)] \rightarrow T$ such that $q(0)\!=\!\sigma $ and $q(\delta(\sigma, \sigma^\prime))=
\sigma^\prime$. In this case, let us write $\llb \sigma, \sigma^\prime\rrb:=q([0, \delta(\sigma, \sigma^\prime)])$, the geodesic from $\sigma$ to $\sigma'$.
\item For all injective continuous functions $g\!: \! [0, 1] \! \rightarrow \! T$, we have $g([0, 1])\! =\! \llb g(0), g(1)\rrb$.
\end{compactenum}
\vspace{0.5mm}
Alternatively, real trees are characterized by the \emph{four-point inequality}:  a connected metric space $(T, \delta)$ 
is a real tree if and only if for all $\sigma_1, \ldots ,\sigma_4\in T$, 
\begin{equation}
\label{4pt}
\delta(\sigma_1, \sigma_2) + \delta(\sigma_3, \sigma_4) \leq \big(\delta(\sigma_1, \sigma_3) + \delta(\sigma_2, \sigma_4)\big) \vee  \big( \delta(\sigma_1, \sigma_4) + \delta(\sigma_2, \sigma_3)  \big) . 
\end{equation}
See Evans \cite{evans05} for more details. In this work, we only consider \emph{compact} real trees.

A \emph{rooted real tree} is a real tree $(T, \delta)$ with a distinguished point $\rho\in T$, called the \emph{root}. 
The \emph{degree} of a point $\sigma \in T$, denoted by $\degree(\sigma)$, is the (possibly infinite) number of connected components of $T\setminus\!\{\sigma \}$. Then, $\sigma$ is said to be a \emph{leaf} if $\degree (\sigma)\! =\! 1$, a \textit{simple point} if $\degree (\sigma)\! =\! 2$ and a \textit{branch point} if $\degree (\sigma) \! \geq\!  3$.

\paragraph{Centers of compact real trees. } Let  $(T, \delta)$ be a compact real tree. We set
\begin{equation}
\label{khvfv}
\forall \, \sigma\in T, \;  \Gamma (\sigma, T)= \sup_{\sigma^\prime \in T} \delta( \sigma, \sigma^\prime ), \quad \text{ then } \quad \diam(T)=\sup_{\sigma \in T} \Gamma (\sigma ,T ),  \quad   \gamma (T) = \inf_{\sigma \in T} \Gamma (\sigma ,T)\; .
\end{equation}
The above definitions make sense as $T$ is compact. If $T$ is rooted at $\sigma$, then $\Gamma(\sigma, T)$ corresponds to the \textit{total height} of the rooted tree $(T, \delta , \sigma)$. In that case, we simply denote  $\Gamma (T)=\Gamma(\sigma, T)$. Note that $ \diam(T) \! =\! \max_{ \sigma, \sigma^\prime \in T}\delta( \sigma, \sigma^\prime )$ is the \emph{diameter} of the metric space $(T, \delta)$. Recall from (\ref{centerdef}) Jordan's definition for the  centers of discrete trees. We introduce an analogous definition for real trees:  
\begin{equation}
\label{def: center}
\text{ a point $c\in T$ is said to be a \textit{center} of the real tree $(T,\delta) $ if } \, \Gamma (c,T) =  \gamma (T).
\end{equation}
We have the following fact about the centers of compact real trees, which is the analog of Jordan's theorem mentioned above. 

\begin{lem} 
\label{realTcentr} Let $(T, \delta)$ be a compact real tree. The following statements hold true:  
\vspace{0.5mm}
\begin{compactenum}[a)]
\item
we have $\gamma (T)\! = \! \frac{1}{2} \diam (T)$\,; 
\item
there exists a unique center $c$ of $(T, \delta)$\,; 
\item
for all pairs of points $(\sigma, \sigma^\prime)$ such that $\delta(\sigma, \sigma^\prime) =  \diam (T)$\,, we have $c \in \llb \sigma, \sigma^\prime \rrb$ and 
$\delta(c, \sigma) = \delta(c, \sigma^\prime) = \frac{1}{2} \diam (T)$\,.
\end{compactenum}
\end{lem}

\noi
\textbf{Proof:} since $T$ is compact, there exist $\sigma, \sigma^\prime \! \in \! T$ such that $\delta(\sigma, \sigma^\prime)\! =\!  \diam (T)$. Let $c\! \in \! \llb \sigma, \sigma^\prime \rrb$ be such that 
$d(c, \sigma)\! =\! d(c, \sigma^\prime)\! =\!  \frac{_1}{^2} \diam (T)$. Let $s \! \in \! T$ be an arbitrary point. We apply the four-point inequality (\ref{4pt}) with $\sigma_1\! =\!  \sigma$, $\sigma_2 \! = \! \sigma^\prime$, $\sigma_3\! =\!  c$ and $\sigma_4\! =\!  s$ and we get 
after some simplification that  
\begin{equation}
\label{djhfou}
\frac{_1}{^2} \diam (T)+ \delta(s, c)\leq \delta(s, \sigma) \vee \delta(s, \sigma^\prime)\le \Gamma(s, T) \; .
\end{equation}
This entails that $\Gamma(s, T)\ge \frac{_1}{^2}  \diam(T)$ for all $s\in T$. It follows that $\gamma(T)\ge \frac{_1}{^2} \diam(T)$. On the other hand, since $\Gamma(s, T)\le \diam(T)$,  \eqref{djhfou} entails that $\delta(c, s)\le \frac{_1}{^2} \diam(T)$ for any $s\in T$. It follows that $\Gamma(c, T)\le \frac{_1}{^2} \diam(T)$. Combined with the previous argument, we obtain that $\gamma(T)=\frac{_1}{^2} \diam(T)=\Gamma(c, T)$. Thus, $c$ is a center of $T$. 
Next, observe that if $s$ is a center of $T$, then $\Gamma(s, T)=\frac{_1}{^2} \diam(T)$. Then we get from \eqref{djhfou} that $\delta(s, c)=0$. 
This then entails the statements in b) and c) and thus completes the proof of the lemma. \cqfd 

\paragraph{The coding of real trees. } 
Let us briefly recall how real trees can be obtained from continuous nonnegative functions. To that end, we write $\bC(\bbR_+, \bbR_+)$ for the space of $\bbR_+$-valued continuous functions equipped with 
the Polish topology of the uniform convergence on every compact subset of $\bbR_+$. Let us denote by $H\! = \! (H_t)_{t\geq 0}$ the canonical process on $\bC(\bbR_+, \bbR_+)$. 
We are concerned with the case where $H$ has a compact support, $H_0 \! =\! 0$ and $H$ is distinct from the null function. We call such a function a \textit{coding function}. Let us assume that $H$ is a coding function. Set $\zeta(H) \! = \! \sup \{ t\! >\! 0: H_t \! >\! 0 \}$, the \textit{lifetime} of the coding function $H$. Note that $\zeta(H) \! \in \! (0, \infty)$ under our assumptions. 
For every $s,t\! \in \! [0, \zeta(H)]$, we set
\[
b_H(s,t)=\inf_{r\in[s\wedge t,s\vee t]}H_r \quad {\rm and} \quad \delta_H(s,t)=H_s+H_t-2b_H(s,t).
\]
It is straightforward to check that $\delta_H$ satisfies the four-point inequality \eqref{4pt}.
Note that $\delta_H$ is a pseudo-metric.  We then introduce the equivalence relation $\sim_H$ by setting 
$s\! \sim_H \! t$ iff $\delta_H(s,t) \! =\! 0$. Let 
\begin{equation}
\label{codef}
\cT_H= [0, \zeta(H)] / \! \sim_H \; . 
\end{equation}
Standard arguments show that $\delta_H$ induces a metric on the quotient set $\cT_H$, which we keep denoting by $\delta_H$. Let  
$p_H \! : \! [0, \zeta(H)] \! \rightarrow \! \cT_H$ be the \textit{canonical projection}. Since $H$ is continuous, so is $p_H$ and $(\cT_H, \delta_H)$ is therefore a compact connected metric space which further satisfies the four-point inequality: it is a compact real tree. We also set $\rho_H\!: =\!  p_H(0) \! =\!  p_H(\zeta(H))$ to be the \textit{root} of $\cT_H$. 
See Duquesne \cite{Duquesne06} for more details on the coding of real trees by functions. 

\paragraph{Re-rooting. } In what follows, we sometimes need to perform a re-rooting operation on the rooted real trees. In terms of the coding functions, this operation corresponds to the following path transformation, which we recall from Duquesne \& Le Gall \cite{Duquesne09}. Let $H$ be a coding function as defined above and recall that its lifetime $\zeta(H) \! \in \! (0, \infty)$. For any $t\! \in \! \bbR_+$, denote by $\overline{t}$ the unique element of $[0, \zeta(H))$ such that $t\! -\! \overline{t}$ is an integer multiple of $\zeta(H)$. If $t_0\in [0, \zeta(H)]$, we define a coding function $H^{[t_0]}$ as follows: 
\begin{equation}
\label{rerootH}
\forall t\! \in \! [0, \zeta(H)] , \quad  H^{[t_0]}_t = \delta_H \big(  \overline{t_0}, \overline{t+t_0}\,  \big) \quad \textrm{and} \quad 
\forall t \geq \zeta(H), \quad   H^{[t_0]}_t =0 \;. 
\end{equation}     
It is not difficult to see that  we can identify $\big( \cT_{H^{[t_0]}}, \delta_{H^{[t_0]}} , \rho_{H^{[t_0]}} \big)$ as the re-rooted tree $\big( \cT_H, \delta_H , p_H(t_0) \big)$. 

\paragraph{Height processes of L\'evy trees and excursion measures. } In \cite{legalllejan}, Le Gall \& Le Jan introduce the height processes, which are the coding functions of L\'evy trees (see also Duquesne \& Le Gall \cite{Duquesne02}). 
We recall here the definition of height processes from their works. We refer the reader to Bertoin's book \cite{Bebook} for background on L\'evy processes. 
Let  $(X_t)_{t\in \bbR_+}$ be a spectrally positive L\'evy process starting from $0$ defined on the probability space $(\Omega, \cF, \bP)$. Then its law is characterized by its Laplace exponent $\Psi\! : \! \bbR_+ \! \rightarrow \! \bbR$ in the sense that
$ \bE [ \exp (-\lam X_t ) ] \! =\! \exp( t\, \Psi (\lam) )$,  for all  $t, \lambda \! \in \! \bbR_+$.  
We assume that $X$ does not drift to $\infty$. In this case, $X_t$ has finite expectation and $\Psi$ takes the following L\'evy-Khintchine form: 
\begin{equation}
\label{LKfor}
 \Psi(\lambda)=a \lambda+b\lambda^2+\int_{(0, \infty)}\pi(dx)\big(e^{-\lambda x}-1+\lambda x\big), \quad \lambda \! \in \!  \bbR_+, 
\end{equation}
where $a, b \! \in \! \bbR_+$ and $\pi$ is a sigma-finite measure on $(0, \infty)$ satisfying $\int_0^\infty (x^2\wedge x)\pi(dx) \! <\! \infty$. 
We restrict our attention to the case  where $\Psi$ satisfies the following assumption: 
\begin{equation}
\label{hypo}
\int_1^\infty \frac{d\lambda}{\Psi (\lam)} < \infty \; .
\end{equation}
Under this assumption, there exists a continuous nonnegative process $H\! =\!  (H_t)_{t\geq 0}$ such that for all $t\! \in \! \bbR_+$, 
the following limit holds in $\bP$-probability: 
\begin{equation}\label{eq: defH}
H_t=\lim_{\ep\to 0}\frac{1}{\ep}\int_0^t ds \,\indi_{\{I^s_t<X_s<I^s_t+\ep\}},
\end{equation}
where $I^s_t:=\inf_{s<r<t}X_r$. The process $H$ is called the \textit{$\Psi$-height process}. In the special case where $\Psi(\lambda)=\lambda^2$ (i.e.\ the Brownian case), classical arguments show that $H$ is distributed as a reflected Brownian motion.

It turns out that $H$ encodes a sequence of random real trees: each excursion of 
$H$ above zero corresponds to a tree in this sequence. 
More precisely, for any $t\in \bbR_+$, we set $I_t= \inf_{s\in [0, t]} X_s$. Basic results of the fluctuation theory entail that $X - I$ is a $\bbR_+$-valued strong Markov process, that $0$ is regular for 
$(0, \infty)$ and recurrent. Moreover, $\! -I$ 
is a local time at $0$ for $X\! -\! I$. 
We denote by $\bN$ the corresponding excursion 
measure of $X\! -\! I$ above $0$.
We can derive from (\ref{eq: defH}) that $H_t$ only depends on the excursion of $X\! -\! I$ above $0$ which straddles $t$. Moreover, we have 
$\{ t\! \in \! \bbR_+ : H_t \! >\! 0 \} \! =\!  \{ t \! \in \! \bbR_+ : X_t \! >\! I_t\}$, that is, the excursion intervals of $H$ above $0$ coincide with those of $X$ above $I$. Let us denote by $(g_i, d_i)$, $i\! \in \! \cI$,  these excursion intervals. Set $H^{i}_s = H_{(g_i +s)\wedge d_i}$, $s\! \in \! \bbR_+$. Then the point measure $\mathcal N:=
 \sum_{i\in \cI} \delta_{(-I_{g_i} ,\,  H^{i})}$   
is a Poisson point measure on $\bbR_+ \! \times \! \bC(\bbR_+, \bbR_+)$ with intensity $dx \, \bN(dH)$. Here,  we have slightly abused the notation by letting $\bN (dH)$ stand for the ``distribution'' of $H(X)$ under $\bN (dX)$, which is a sigma-finite measure on $\bC(\bbR_+, \bbR_+)$.   
In the Brownian case, up to a multiplicative constant, $\bN$ is the Ito's positive excursion measure of  Brownian motion and $\mathcal N$ reduces to the Poisson decomposition of a reflected 
Brownian motion above $0$. 

In the rest of the paper, we will work exclusively with the $\Psi$-height process 
$H$ under the excursion measure $\bN$. 
The following holds true:
\[
 \textrm{$\bN$-a.e.} \quad \zeta(H)  <  \infty\, , \quad H_0 = H_{\zeta(H)} = 0 \quad \textrm{and} \quad H_t  > 0 \;  \Longleftrightarrow \; t  \in  (0, \zeta(H)) \; .
\]
This shows that $H$ under $\bN$ is a coding function as defined above. Duquesne \& Le Gall \cite{Duquesne05} then define the \textit{$\Psi$-L\'evy tree} as the real tree coded by $H$ under $\bN$ in the sense of \eqref{codef}. 
In that case, when there is no risk of confusion, we simply write $ \big( \cT, \delta, \rho \big)$ instead of  $( \cT_H, \delta_H, \rho_H )$. 

\paragraph{L\'evy trees conditioned by their total heights. } Let $H$ be the $\Psi$-height process under its excursion measure $\bN$, as defined above. 
We use the notation $\Gamma(H):= \sup_{t \in [0, \zeta(H)]} H_t$, which coincides with the total height of the $\Psi$-L\'evy tree $(\cT, \delta, \rho)$.  
Let us recall from Duquesne \& Le Gall \cite{Duquesne02} (Corollary 1.4.2) the following distributional properties of $\Gamma(H)$. 
\begin{equation}
\label{dist_Gamma}
\forall \,r \! \in \! (0, \infty), \quad \bN\big(\,\Gamma(H)  >   r\big)= v(r)  \, , \quad  \textrm{where $v$ verifies} \; \int_{v(r)}^\infty \frac{d\lambda}{\Psi (\lam)} = r\; ; 
\end{equation}
also $\bN$-a.e., there exists a unique $t \! \in \! (0, \zeta(H)) $ such that $H_t=\Gamma(H)$.

Abraham \& Delmas in \cite{AbDe09} define the laws of the $\Psi$-height process conditioned by the total heights. More precisely, they construct a family of probability laws $\bN (\,\cdot\,| \, \Gamma(H) \! =\!  r)$, $r\! \in \! (0, \infty)$, on $\bC(\bbR_+, \bbR_+)$ which satisfy the following properties: 
\begin{compactenum}[-]
\item[(a)]
$\bN (\,\cdot\,| \, \Gamma(H) \! =\!  r)$-a.s.~$\Gamma(H) \! = \! r$. 
\item[(b)] The mapping $r\mapsto \bN (\,\cdot\,| \, \Gamma(H) \! =\!  r)$ is continuous with respect to the weak topology on $\bC(\bbR_+, \bbR_+)$. 
\item[(c)] $\bN = \int_0^\infty  \bN(\Gamma(H)   \in dr)\,  \bN (\,\cdot\,| \, \Gamma(H) \! =\!  r)$. 
\end{compactenum}
In addition, the authors of \cite{AbDe09} give a Poisson decomposition at the unique maximum point of $H$, which generalizes Williams' decomposition for Brownian excursions.

\paragraph{L\'evy trees conditioned by their diameters. } Recall that $(\cT, \delta, \rho)$ stands for the L\'evy tree coded by the $\Psi$-height process $H$ under its excursion measure $\bN$. In~\cite{DuWa14}, 
Duquesne \& W.\ study the diameter of $\mathcal T$ as well as a spinal decomposition of $\mathcal T$ along its longest geodesic. In particular, the following results can be found there. 
We have that $\bN$-a.e.~there exists a unique pair of times $0\! < \! \tau_0\! < \! \tau_1 \! <\! \zeta(H) $ such that 
$\delta (\tau_0, \tau_1)=\diam(\cT)$. The distribution of $\mathrm D=\mathrm D(H)\!:=\diam(\cT)$ under $\bN$ is given by
\[
\bN \big( \mathrm{D} \! > 2r  \big)= v(r) -\Psi \big( v(r)\big)^2 \!\! \int_{v(r)}^\infty \! \frac{d\lambda}{\Psi(\lambda)^2}\; , \quad r\in (0, \infty) \; ,
\]
where $v$ is defined in \eqref{dist_Gamma}. 
Next, we introduce the laws of $\cT$ conditioned by its diameter. To that end, let $H, \widetilde{H} \! \in \! \bC(\bbR_+, \bbR_+)$ be two coding functions. 
The \textit{concatenation} $H\oplus \widetilde H$ of $H$ and $\widetilde{H}$ is the coding function defined as 
\begin{equation}
\label{concadef}
\forall t\in \bbR_+, \qquad (H\oplus \widetilde{H})_t = H_t \quad \textrm{if $t\in [0, \zeta(H)]$} \quad \textrm{and} \quad (H\oplus \widetilde{H})_t = \widetilde{H}_{t-\zeta(H)} \quad \textrm{if $t \geq \zeta(H)$.}
\end{equation}
For all $r \! \in \! (0, \infty)$, we define $\bQ_r$ as the law of $H  \oplus  \widetilde{H}$ under $\bN (dH\,| \, \Gamma(H) \! =\!  r/2)\bN (d\widetilde{H}\,| \, \Gamma(\widetilde{H}) \! =\!  r/2)$. 
Namely, for all measurable functions $F \! :\!  \bC (\bbR_+, \bbR_+) \! \rightarrow \! \bbR_+$,
\begin{equation}
\label{defQr}
\bQ_r \big[ F(H)\big]= \int \!\!\!\! \int_{\bC(\bbR_+, \bbR_+)^2}\!\!\!  F\big( H \!  \oplus \! \widetilde{H} \big)\; \bN\big( dH \, \big| \, \Gamma(H) =  r/2\big) \, \bN\big( d\widetilde{H} \, \big| \, \Gamma(\widetilde{H})  =  r/2\big)\;   \; .
\end{equation}
Then $\bQ_r$ has the following properties.  
\begin{compactenum}
\item[(a)] $\bQ_r$-a.s.~we have $\mathrm{D}\! =\! r$ and there exists a unique pair of points $\tau_0, \tau_1\! \in \! [0, \zeta(H)]$ such that $\mathrm{D} \! = \! \delta (\tau_0, \tau_1)$. Moreover, the unique center of $\cT$ 
has degree $2$: it is a simple point.   
\item[(b)] For all $r\! \in \! (0, \infty)$, $\bQ_r  [ \, \zeta(H) \, ]  \! \in \! (0, \infty)$. Moreover, the mapping  
$r\! \mapsto \! \bQ_r$ is weakly continuous and for all measurable functions $F\! :\!  \bC (\bbR_+, \bbR_+) \! \rightarrow\!  \bbR_+$ and $f \! :\!  \bbR_+ \! \rightarrow \! \bbR_+ $, 
\begin{equation}
\label{desindiam}
 \bN \big[ f(\mathrm{D}) F(H)\big]= \int_{0}^\infty  \! \frac{\bN (\,\mathrm{D} \in  dr)}{\bQ_r[\, \zeta(H)\, ]}
\, f(r) \,  \bQ_r \Big[ \int_0^{\zeta(H)} \! F\big(  H^{[t]} \big) \, dt \Big]    \; , 
\end{equation} 
where $H^{[t]}$ the rerooting of $H$ at $t$, defined in (\ref{rerootH}). 
\end{compactenum}
It follows from \eqref{desindiam} that we have a regular version of the conditional laws $\bN (dH  \,|\, \mathrm{D}\! = r\! \, )$ that are obtained from $\bQ_r$ by a uniform re-rooting:  
for all measurable functions $F \! :\!  \bC(\bbR_+, \bbR_+) \! \rightarrow \! \bbR_+$, we have 
\begin{equation}
\label{Hcondia}
\forall \,r \! \in \! (0, \infty), \qquad \bN \big[ \,F(H) \, \big| \, \mathrm{D} \! = r\! \, \big]  =   \bQ_r \Big[ \int_0^{\zeta(H)} \! F\big(  H^{[t]} \big) \, dt \Big]  \Big/ \bQ_r[\, \zeta(H)\, ] \; . 
\end{equation}
We call $\bQ_r$ the law of \textit{unrooted $\Psi$-L\'evy trees with diameter $r$}. For a more extensive account, we 
refer to Duquesne \& W.~\cite{DuWa14}, Theorem 1.2 and Remark 1.2.   
 
\section{Main results } 
\label{sec: result}

Recall that all the random variables here are defined on the same probability space $(\Omega, \cF, \bP)$. 
In particular, it contains the following. 
\begin{compactenum}[--]
\item
A spectrally positive L\'evy process $(X_s)_{s \in \bbR_+}$ whose Laplace exponent $\Psi$ is given in \eqref{LKfor} and satisfies the condition (\ref{hypo}). 
\item
For all positive integers $p$, let $\mu_p\! = \! (\mu_p (k))_{k\in \N}$ be a probability 
law on $\N=\{0, 1, 2, \dots\}$ that verifies (\ref{subcrGW}). Let $(J^{(p)}_n)_{n\geq 1}$ be a sequence of independent random variables with common law $\mu_p$. 
We assume that there exists a non decreasing sequence of positive integers $(b_p)_{p\geq 1}$ such that the following convergence holds in distribution for $\bbR$-valued random variables: 
\begin{equation}
\label{assH}
\frac{p}{b_p} \big( J^{(p)}_1 + \ldots + J^{(p)}_{b_p} -b_p\big) \overset{\textrm{(d)}}{\underset{p\rightarrow \infty}{-\!\!\!- \!\!\! \longrightarrow}} X_1  \; .  
\end{equation}
\end{compactenum}
We set $ g^{\mu_p} (s)  \! =\!  \sum_{k\in \N} s^k \mu_p (k)$, $s\! \in \! [0, 1] $, the generating function of $\mu_p$.  Let us denote by $g^{\mu_p}_n$ the $n$-th iteration of $g^{\mu_p}$, that is, 
$g^{\mu_p}_{n+1} \! = \! g^{\mu_p}  \circ  g^{\mu_p}_n \! =\!  g^{\mu_p}_n \circ g^{\mu_p}$ and $g^{\mu_p}_0 (s)\! = \! s$, $s\! \in \! [0, 1]$. We assume that for all $r\! \in \! (0, \infty) $,  
\begin{equation}
\label{asszeta}   
\liminf_{p\rightarrow \infty} \big( g^{\mu_p}_{\lfloor p r \rfloor } (0) \big)^{b_p/p} > 0 \; .
\end{equation}
Let $\tau_p \! : \! \Omega \rightarrow  \Tor$ be a GW($\mu_p$)-tree. 
Recall from Section \ref{sec: pl} the contour function $(C_s (\tau_p), 0\le s\le 2(|\tau_p| -1))$ of the ordered rooted tree $\tau_p$. 
For convenience, we extend the definition of $C(\tau_p)$ to $\bbR_+$ by setting $C_s(\tau_p)\! =\!  0$ for all $s \! \geq \! 2(|\tau_p| \! -\! 1)$. 
We also set 
\[
\Gamma (\tau_p)  =  \sup_{0\le s\le 2(|\tau_p| -1)} C_s (\tau_p)\,,
\] 
which coincides with the total height of $\tau_p$\,. Recall from Section \ref{sec: Levy} the $\Psi$-height process $H$ defined under the excursion measure $\bN$. 
Under the assumptions (\ref{assH}) and (\ref{asszeta}), Duquesne \& Le Gall in~\cite{Duquesne02} have shown a general invariance principle for $(\tau_p)_{p\ge 1}$\,; see Theorem 2.3.1 and Corollary 2.5.1 there. Then, they deduce (Proposition 2.5.2 of \cite{Duquesne02}) that  for all $r\in (0, \infty)$,  
\begin{equation}
\label{heigthcjkb} 
\big(   p^{-1} C_{2 b_p s} (\tau_p), s\in \bbR_+\big) \; \; \text{under} \; \; \bP (\, \cdot \, | \, \Gamma (\tau_p)  \geq   pr) \quad \xrightarrow[p\to\infty]{\text{(d)}} \quad  H \; \; 
\text{under} \; \; \bN (\, \cdot \, | \, \Gamma(H) \geq  r) , 
\end{equation}
where the convergence holds in distribution on $\bC (\bbR_+, \bbR_+)$. 
Here, we show that their result can be extended to the following. 

\begin{prop}
\label{bdkjfvs} Let $\Psi$ be given in (\ref{LKfor}) and satisfy (\ref{hypo}). Let $\mu_p$ be a probability law on $\mathbb N$ which verifies \eqref{subcrGW}, for each $p\ge 1$. 
Suppose that (\ref{assH}) and (\ref{asszeta}) take place. 
Let $(b_p)_{p\geq 1}$ be as in (\ref{assH}) and let $r \! \in \! (0, \infty)$. 
Let $C(\tau_p)$ be the extended contour function of the GW$(\mu_p)$-tree $\tau_p$ as defined above. 
Then, the following convergence holds in distribution on $\bC (\bbR_+, \bbR_+)\times \bbR_+$:
\begin{equation}
\label{heightcon} 
\Big(\!\big(p^{-1}  C_{2 b_p s} (\tau_p), s\in \bbR_+\big)  ,  \tfrac{|\tau_p|}{b_p}  \Big)\; \textrm{under} \; \bP \big(\, \cdot \, \big| \, \Gamma (\tau_p) \! =\!  \lfloor pr\rfloor \big)\;  \overset{{\textrm{(d)}}}{\underset{p\rightarrow \infty}{-\!\!\!-\!\!\!-\!\!\!- \!\!\! \longrightarrow}} \;  \big( H , \zeta(H)\big)  \; 
\textrm{under} \; \bN(\, \cdot\,| \, \Gamma(H) \! =\!  r) . 
\end{equation}
Here, $\bN (\,\cdot\, | \, \Gamma(H) \! = \! r)$ stands for the law of the $\Psi$-height process $H$ conditioned to have a total height $r$ and $\zeta(H)$ stands for its lifetime. 
\end{prop}
\begin{rem}
\label{nonsingre}
The assumptions (\ref{assH}) and (\ref{asszeta}) are minimal for (\ref{heigthcjkb}) to hold: see the discussion right after Theorem 2.3.1 in \cite{Duquesne02}, p.~55. 
Let us also mention that if $\mu_p=\mu$ for all $p$ and \eqref{assH} hold, then \eqref{asszeta} is automatically verified. In this case, the limit $X_1$ in \eqref{assH} is necessarily distributed as a spectrally positive $\alpha$-stable random variable, for some $\alpha\in (1, 2]$. See Theorem 2.3.2 in \cite{Duquesne02} for the details. 
\cq
\end{rem}

\bigskip

The aim of this work is to prove a limit theorem for a family of random unrooted unlabelled plane trees which are defined as follows.  
Let $\mu \! = \! (\mu(k))_{k\in \N}$ be a probability law on $\mathbb N$ which satisfies (\ref{subcrGW}). For a plane tree $t\in \bbT_{\! {\rm pl}}$, we set
\begin{equation}
\label{jpkgkjh}
W_{\! \mu} (t)=  \prod_{v \; \textrm{vertex of } t}\!\! \mu \big( \degree (v)\! - \! 1 \big)\, .
\end{equation}
Recall that $D(t)$ stands for the diameter of a tree $t$. For all $k\in \mathbb N$, let us denote
\begin{equation}
\label{def: Tplk}
\Tpl (k)= \big\{ t  \in  \Tpl: \diam (t) =  k \big\} \quad \textrm{and} \quad Z_k (\mu) = \!\!
  \sum_{t\in \Tpl(k)} \! \! W_{\! \mu} (t) \; .
\end{equation}  
Note that (\ref{subcrGW}) implies $Z_k (\mu) \! >\! 0$ for all $k\in \mathbb N$. Though 
$ \bbT_{\! {\rm pl}} (k)$ has infinite cardinality for each $k\ge 2$, we prove in Lemma \ref{Zpfini} that $Z_k (\mu) \! <\! \infty $. As a result, the following probability law is well defined for each $k$ : 
\begin{equation}
\label{fljgrd}
\forall t \in   \Tpl (k), \quad  Q^\mu_k (t) = \frac{W_{\! \mu} (t)}{Z_k(\mu)} \; .
\end{equation}

The above Proposition \ref{bdkjfvs} plays an important role in our study of  $Q^\mu_k$. Indeed, by rooting plane trees at their central edges (see the definition below), we will see (in Lemma \ref{contigu}) that $Q^\mu_k$ is closely related to the laws of Galton--Watson trees conditioned by total heights. 

\paragraph{Central edges.} \label{p: ed}Let $t\in\Tpl$ be a plane tree with vertex set $V$. Recall from \eqref{centerdef} Jordan's definition for the centers of a tree. 
\begin{compactenum}[--]
\item If $\diam(t)$ is odd, then $t$ has two adjacent centers $c, c^\prime$. In this case, we define the \textit{set of central edges of $t$} as $K(t)= \{ (c,c^\prime), (c^\prime, c) \}$. Namely, $t$ has exactly two central edges, which are the two oriented edges between its two centers. 
\item If $\diam(t)$ is even, then $t$ has a unique center $c$. In this case, we define the \textit{set of central edges of $t$} as 
$K(t) \! = \! \big\{ ( v, c)\, ; \,  v \! \in \! V : v \! \sim \! c \; \textrm{and} \; \exists\, w \! \in \! V:  \diist(w, c) \! =\!  \frac{_1}{^2}\diam (t) \; \textrm{and} \; v \! \in \! \llb w, c \rrb  \big\}$. 
In other words, a central edge is an oriented edge $(v, c)$ of $t$ where $v$ belongs to a path of length $D(t)$.  
Clearly, we have $2\! \leq \! |K(t)| \! \leq \! \degree (c)$. 
\end{compactenum} 
In particular, observe that the head of a central edge is necessarily a center of $t$. 

\medskip
Let $r\in (0, \infty)$. Let $(\mu_p)_{p\ge 1}$ be  a sequence of probability measures on $\mathbb N$ that verify \eqref{subcrGW}. For each $p$, let $\Trp$ be a random plane tree  whose distribution is given by $Q^{\mu_p}_{\lfloor rp\rfloor}$, as defined in (\ref{fljgrd}), that is, 
\[
\bP \big(\Trp  =  t\big)  =  Q^{\mu_p}_{\lfloor rp \rfloor} (t)=\frac{1}{Z_{\lfloor rp \rfloor}(\mu_p)}\prod_{v \; \textrm{vertex of } t}\mu_p\big( \degree (v)\! - \! 1 \big), \quad t \in  \Tpl (\lfloor rp \rfloor)\; .
\]
Given $\Trp$, let $\epp$ be a central edge picked uniformly from $K(\Trp)$. We root $\Trp$ at $\epp$, giving rise to an edge-rooted plane tree $(\Trp\,, \epp)$. Recall from Section \ref{sec: pl} its contour function $C(\Trp\,, \epp)$ defined on $[0, 2(|\Trp|-1)]$. We extend its definition by setting $C_s(\Trp,\, \epp)=0$ for $s\ge 2(|\Trp|-1)$. 
 
\begin{thm}
\label{main} 
Let $\Psi$ be given in (\ref{LKfor}) and satisfy (\ref{hypo}). Let $(\mu_p)_{p\ge 1}$ be a sequence of probability laws that verify \eqref{subcrGW}. 
Suppose that  (\ref{assH}) and (\ref{asszeta}) take place.   
Let $(b_p)_{p\geq 1}$ be as in (\ref{assH}). For $r\in (0, \infty)$, let $C(\Trp, \epp)$ be the extended contour function of the edge-rooted plane tree $(\Trp, \epp)$ as defined above. Then 
the following convergence holds in distribution on $\bC (\bbR_+, \bbR_+)$:
\begin{equation}
\label{heigthcon} 
\big( p^{-1} C_{2 b_p s} (\Trp\,, \epp),\, s\in \bbR_+\big) \ \textrm{ under } \; \,  \bP \quad  \xrightarrow[p\to\infty]{\text{(d)}} \quad   H \; 
\textrm{ under } \; \bQ_r \,, 
\end{equation}
where $\bQ_r$ stands for the law of unrooted $\Psi$-L\'evy trees with diameter $r$ as defined in (\ref{defQr}). 
\end{thm}

\begin{rem}
\label{rem: main}
By standard arguments (see for instance \cite{AbDeHo13}), the convergence in \eqref{heigthcon} implies the weak convergence of $(\Trp, p\ge 1)$ in Gromov--Hausdorff--Prokhorov topology. Indeed, for each $p\ge 1$, denote by $(\Trp, p^{-1}\diist)$ the metric space obtained from the graph $\Trp$ after rescaling its graph distance $\diist$ by $1/p$. Let $\frac{1}{b_p}\mathrm m^{_{(p)}}$ be the (finite) measure of $\Trp$ obtained by putting a mass $b_p^{-1}$ at each node of $\Trp$. 
Recall from Section \ref{sec: Levy} the real tree $(\cT_H, \delta_H)$ encoded by the canonical process $H$ on $\bC(\bbR_+, \bbR_+)$. Denote by $\mathrm m_H$ the push-forward measure of the Lebesgue measure on $[0, \zeta(H)]$ by the projection $p_H: [0, \zeta(H)]\to \cT_H$. 
Then, under the assumptions of Theorem \ref{main}, we have 
\[
\Big(\Trp, \tfrac{1}{p}\,\diist, \tfrac{1}{b_p}\,\mathrm m^{_{(p)}} \Big)\xrightarrow[p\to\infty]{(d)} \big(\mathcal T_H, \delta_H, \mathrm m_H\big) \;  \text{ under } \;\bQ_r,
\]
with respect to the Gromov--Hausdorff--Prokhorov topology. 
Similarly, we can reformulate the convergence in \eqref{heightcon}  in terms of a Gromov--Hausdorff--Prokhorov  convergence of the conditioned Galton--Watson trees. 
\end{rem}

The rest of the paper is organized in the following way. In
Section \ref{Pfbdkjfvs}, we provide the proof of Proposition \ref{bdkjfvs}, based on the convergence in \eqref{heigthcjkb} and the following observation: 
take a Galton--Watson tree conditioned to have a total height at least $p$ and locate its first tip (i.e.~the first node at maximum height in lexicographic order); step down along the ancestral line of this tip to a depth $p$; taking the path of length $p$ along with all the trees planted on it gives a subtree of the initial tree; it turns out that this subtree is distributed as a Galton--Watson tree conditioned to have a total height $p$. 
See (\ref{defthetre}) and Lemma \ref{tiplaw} for a precise statement. Section \ref{mainthpf} is devoted to the proof of Theorem \ref{main}, where we employ the following idea. 
For the unrooted L\'evy tree with diameter $r$, a decomposition at its unique center 
yields two independent (rooted) L\'evy trees with total height $r/2$ (this can be considered as a verbal description of the definition \eqref{defQr}; see also the point (a) right below it). 
Then the main point of the proof is to show that asymptotically as $p\to \infty$, the random plane trees $\Trp$ also behaviors in a similar way, which is achieved by establishing an upper bound for the number of its central symmetries.

\section{Proof of Proposition \ref{bdkjfvs}}
\label{Pfbdkjfvs}
Recall the notation $\bbU$ from Section \ref{sec: or}. 
If $u\! = \! (a_1, \ldots, a_n) \in \bbU$, we write $|u|= n$ for the \emph{length} of $u$. Let $v= (b_1, \ldots , b_m) \in \bbU$. We denote by $u\ast v= (a_1, \ldots, a_n , b_1, \ldots , b_m) $ the \textit{concatenation} of $u$ with $v$. 
Let $t\in \bbT_{\! {\rm or}}$ be a finite ordered rooted tree, which is a subset of $\bbU$. 
For $u  \in  t$, we define the \textit{subtree of $t$ stemming from $u$} as 
\begin{equation}
\label{thetastem}
 \theta_u (t) \! = \! \big\{ v \! \in \! \bbU: u\ast v  \in t \big\} \; .
 \end{equation}
Observe that $\theta_u (t)  \in   \Tor$. 

The set $\bbU$ is naturally equipped with the \textit{lexicographical order} denoted by $\preceq$, which is a total order on $\bbU$. Since $t$ is finite, $\preceq$ induces a linear order on it.   
Recall that 
$\Gamma (t)\! = \! \max_{v\in t} |v|$ stands for the total height of the tree $t$.
We say that a vertex $v\in t$ is a \textit{tip of $t$} if  
$|v|\! =\!  \Gamma (t)$. 
Let $u$ be the 
tip of $t$ which is minimal with respect to the order $\preceq$. 
Suppose that $n, p$ are two integers such that $\Gamma(t)=n\ge p\ge 0$.
Then we can write $u=(a_1, a_2, \dots, a_n)$ for some $a_j\in \mathbb N^*$, $1\le j\le n$. 
Let us denote by $u_0=\varnothing$ and by $u_j=(a_1, a_2, \dots, a_j)$ for $1\le j\le n$. In other words, the sequence $\{ u_0, \ldots, u_n \}$ forms the ancestral line of $u$. 
We define 
\begin{equation}
\label{defthetre}
\Theta_p (t)= \theta_{u_{n-p}} (t)\, , \quad U_p (t)= u_{n-p} \, , \quad \textrm{and} \quad \Lambda_p (t)= \{ u_{n-p} \} \cup \big( t \backslash  \Theta_p (t) \big) \; .
\end{equation}
Note that both $\Theta_p (t)$ and  $\Lambda_p (t)$ are ordered rooted trees and that $\Gamma(\Theta_p (t))=p$. We have the following.
 
\begin{lem}
\label{tiplaw} Let $\mu\! = \! (\mu(k))_{k\geq 0}$ be a probability distribution satisfying (\ref{subcrGW}). Let $\tau$ be a GW($\mu$)-tree defined on the probability space $(\Omega, \mathcal F, \bP)$. Then for any $p\in \mathbb N$, we have  
$$\tau \ \textrm{ under } \ \bP ( \, \cdot \, | \, \Gamma (\tau)  = p )  \quad \eqd \quad   \Theta_p (\tau) \ \textrm{ under } \ \bP ( \, \cdot \, | \, \Gamma (\tau)  \geq  p ) \; .$$
\end{lem}
\noi
\textbf{Proof:} we set respectively $\Tor^{^=} (p)  =  \{ t \in  \Tor : \Gamma (t) =  p \}$, $\Tor^{^\geq } (p)  =  \{ t \in  \Tor : \Gamma (t) \geq   p \}$ and 
\[
A_p \! =\!  \big\{ (t,u) ; t  \in   \Tor , u  \in  t :  k_u(t)=0 \, \text{ and } \, \forall\, v \in  t,\,   |v|  <  |u|+ p+ \un_{\{ u\preceq v\}}   \big\} \; ,
\]
where we recall the notation $k_u(t)$ standing for the number of children of $u$ in the family tree $t$. 
We observe that  
\[
t   \in   \Tor^{^\geq } (p) \longmapsto  \Big( \Theta_p (t) \, , \, \big( \Lambda_p (t), U_p (t)\big) \Big) \in  \Tor^{^=} (p)  \times  A_p
\]
is a bijective mapping. We denote by $\phi$ its inverse. 
Let $t_1 \! \in \!  \Tor^{^=} (p)$ and let $(t_2, u) \!  \in \! A_p$\,. We deduce from \eqref{GWlawdef} that 
\[ 
\bP \big( \tau =  \phi \big(t_1, (t_2, u) \big) \big) = \bP (\tau =  t_1) \bP (\tau  =  t_2 ) / \mu (0) \; .
\]
Thus,  $\bP \big( \Theta_p (\tau) \! = \! t_1\, ; \, \Gamma (\tau) \! \geq \! p \big)\! = \!  \bP (\tau\! = \! t_1) S_p$, where $S_p= \sum_{(t, u)\in A_p} \bP (\tau \! =\!  t ) / \mu (0)$. 
Summing over all $t_1\in \Tor^{^=} (p)$, we find that $\bP (\Gamma(\tau)\ge p)=S_p\sum_{t_1\in \Tor^{^=} (p)}\bP(\tau=t_1)$. Then the desired result readily follows. 
\cqfd

\bigskip

Next, we reformulate the mapping $t\mapsto \Theta_p (t)$ as a transform of contour functions. To that end,  recall that $\bC(\bbR_+, \bbR_+)$ stands for the set of  continuous functions from $\bbR_+$ to $\bbR_+$. We denote by $\bC$ the set of  coding functions, namely the set $\bC$ consists of  the functions $H\! \in \! \bC (\bbR_+, \bbR_+)$ with compact supports, not identically null and satisfying $H_0\! = \! 0$. Let $H\in\bC$. Then $\Gamma (H) = \sup H \in (0, \infty)$. We set 
\[
 S(H) :=   \inf \big\{ t \in  \bbR_+: H_t  =  \Gamma (H) \big\} \in (0, \zeta(H)),
\] 
where we recall the notation $\zeta(H)=\sup \{ t\in \bbR_+ : H_t \! >\! 0 \}\in (0, \infty)$. 
For all $r  \in  (0, \Gamma (H) )$, the following quantities are well defined: 
\[ 
\sigma^-_r (H) =  \sup\! \big\{ t  \!\in\!  [0, S(H)]: H_t \! <\! \Gamma (H) \! -\! r \big\} \quad \textrm{and}  \quad  
\sigma^+_r (H) =  \inf\! \big\{ t  \!\in\!  [S(H), \infty): H_t \! <\! \Gamma (H) \! -\! r \big\}\,.
\]
We set for all $H\in \bC$ and $r\in (0, \Gamma(H))$\,, 
\begin{equation}
\label{defthefun}
 \Theta_r (H) (t) : = H \big( (\sigma^-_r (H) + t)\wedge \sigma^+_r(H) \big) -\Gamma (H) + r\,, \quad t\in\mathbb R_+\,.
\end{equation} 
Clearly, $\Theta_r (H) \in \bC$. Now let $(\cT_H, \delta_H, \rho_H)$ be the rooted real tree coded by $H$ as explained in \eqref{codef}. Set $\sigma \! :=\!  p_H (S(H))$, the first tip of $\cT_H$, and let $\sigma_r$ be the unique point of the 
geodesic $\llb \rho_H, \sigma\rrb $ satisfying $\delta_H (\sigma_r, \sigma)\! = \! r$. We set $\Theta_r (\cT_H)\! := \! \{ s \! \in \! \cT_H : \sigma_r \! \in \! \llb \rho_H , s \rrb\}$, the subtree of $\cT_H$ stemming from $\sigma_r$.  Then
the tree coded by $\Theta_r (H)$ is isometric to the rooted compact real tree $(\Theta_r (\cT_H), \delta_H, \sigma_r)$. 

In the case of (discrete) ordered rooted tree, we have a similar observation.  Indeed, let $t \! \in \! \bbT_{\! {\rm or}}$ and let $p$ be a positive integer such that $\Gamma (t) \! \geq \! p$. Recall the contour function $(C_s(t))_{s\in \bbR_+}$ of $t$. Then, 
\begin{equation}
\label{contoucsq}
\Theta_p \big((C_s(t))_{s\in \bbR_+} \big) =  \big( C_s(\Theta_p(t) ) \big)_{s\in \bbR_+} \; .
\end{equation}
Here, the first $\Theta_p$ is defined in (\ref{defthefun}) and the second one in (\ref{defthetre}). 

We will need some continuity properties of the mapping $(r, H)\mapsto \Theta_r (H)$. Let us recall that $\bC (\bbR_+, \bbR_+)$ is equipped with the Polish topology induced by the uniform convergence on every compact subset. 
\begin{lem}
\label{conThefun} Let $H, H^{(p)} \!  \in\!  \bC $, $p\! \in \! \N$. Let $r, r_p \! \in \! (0, \infty)$ be such that $r \! < \! \Gamma (H)$ and $r_p \! < \! \Gamma (H^{(p)})$ for each $p$. 
We assume that  the following conditions hold true. 
\begin{itemize}
\item[(i)] $ \{ s \! \in \! \bbR_+ : H_s\! = \! \Gamma (H) \}= \{ S(H) \} $.   
\item[(ii)] For all $s\! \in \! ( \sigma^-_r (H), \sigma^+_r (H))$\,, $H_s  > \Gamma (H)  - r$. 
\item[(iii)] $\lim_{p \rightarrow \infty} H^{(p)} =  H$ in $\bC(\bbR_+, \bbR_+)$\,, $\lim_{p\rightarrow \infty} \zeta(H^{(p)}) =  \zeta(H)$ and $\lim_{p\rightarrow \infty} r_p =  r$\,. 
\end{itemize}
Then,  
\begin{equation}
\label{cvtemps}
\lim_{p \rightarrow \infty} S(H^{(p)}) =  S(H) \, , \quad \lim_{p \rightarrow \infty} \sigma^-_{r_p}(H^{(p)}) =  \sigma^-_{r}(H) \quad \textrm{and} \quad  \lim_{p \rightarrow \infty} \sigma^+_{r_p}(H^{(p)}) =  \sigma^+_{r}(H)\; .
\end{equation}
Moreover, we have 
\begin{equation}
\label{eq: cvTheta}
\lim_{p\rightarrow \infty} \Theta_{r_p} \big(H^{(p)}\big)  =  \Theta_r (H)  \ \text{ in } \bC(\bbR_+, \bbR_+)  \quad \text{ and }\quad \zeta\big(\Theta_{r_p} (H^{(p)})\big)\to \zeta\big(\Theta_r(H)\big)\,. 
\end{equation}
\end{lem}
\noi
\textbf{Proof:} by $(iii)$, there exists some $a\! \in \! (0, \infty)$ such that $ 1+ \zeta(H) +\sup_{p\geq 0} \zeta(H^{(p)})  \leq  a$. Thus, 
$\Gamma (H) =  \max_{[0, a]} H  =  
\lim_{p\rightarrow \infty} \max_{[0, a]} H^{(p)}   =  \lim_{p\rightarrow \infty} \Gamma (H^{(p)})$. 
Let $\ep \! \in \! (0, S(H))$. By $(i)$, we obtain that 
\[
\sup  \big\{ H_s \,; s \! \in \! [0, a] : |s\! -\! S(H)|  \geq  \ep \big\}  <  \Gamma (H)\; .
\] 
Therefore, for all sufficiently large $p$, we have $\sup \{ H^{(p)}_s ; s \! \in \! [0, a] : |s\! -\! S(H)| \! \geq \! \ep \} <  \Gamma (H^{(p)})$. 
Since  $\zeta(H^{(p)})  \leq  a$, this entails that $|S(H^{(p)})\! -\! S(H)| \! \leq \! \ep$. Since $\ep$ can be arbitrarily small, we get $\lim_{p \rightarrow \infty} S(H^{(p)})\! = \! S(H)$. This proves the first convergence in \eqref{cvtemps}.
 
Let $\ep \! \in \! (0, \sigma^-_r (H))$. By definition, there exist $s_1 \! \in \! [ \sigma^-_r (H) \! -\! \ep, \sigma^-_r (H)]$ and $s_2 \! \in \! [ \sigma^+_r (H), \sigma^+_r (H) + \ep]$ such that 
$ \max (H_{s_1} , H_{s_2}) \!  < \!  \Gamma (H) \! - \! r$. Since $\sigma^-_r (H) \! < \! S(H) \! < \! \sigma^+_r (H)$ and $\lim_{p \rightarrow \infty} S(H^{(p)})\! = \! S(H)$, the following holds true for all sufficiently large $p$:
\[
s_1 < S(H^{(p)}) \, , \quad H^{(p)}_{s_1}  <  \Gamma (H^{(p)})  -  r_p \, ,  \quad s_2 > S(H^{(p)} )\quad \textrm{and}  \quad H^{(p)}_{s_2}  <  \Gamma (H^{(p)})  - r_p \, , 
\]
which implies that $  s_1  <   \sigma^-_r (H^{(p)})$ and $ s_2  >   \sigma^+_r (H^{(p)})$. As $\ep$ can be chosen arbitrarily small, we see that
\begin{equation}
\label{facidhf}
\liminf_{p\rightarrow \infty} \sigma^-_r (H^{(p)}) \geq \sigma^-_r (H) \quad \textrm{and} \quad \limsup_{p\rightarrow \infty} \sigma^+_r (H^{(p)}) \leq \sigma^+_r (H) \; .
\end{equation}
Let $0 \!< \! \ep \! < \! \min ( S(H)\! -\! \sigma^-_r (H) , \sigma^+_r (H)\!- \! S(H) ) $; then $S(H) \! \in \!  ( \sigma^-_r (H) \! +\! \ep, \sigma^+_r (H) \! -\! \ep)$.  By (ii), we have 
$ \min\{ H_s\, ; s \in  [ \sigma^-_r (H) \! +\! \ep, \sigma^+_r (H) \! -\! \ep] \}  > \Gamma (H) \! -\! r$. By the fact that  $\lim_{p \rightarrow \infty} S(H^{(p)})\! = \! S(H)$ and by $(iii)$, we deduce that  for all sufficiently large $p$, 
$$  \sigma^-_r (H) \! +\! \ep < S(H^{(p)} ) < \sigma^+_r (H) \! -\! \ep  \quad \textrm{and} \quad  \min\{ H^{(p)}_s ; s \in  [ \sigma^-_r (H) \! +\! \ep, \sigma^+_r (H) \! -\! \ep] \} \! >\! \Gamma (H^{(p)}) \! -\! r_p \; .$$
It follows that $\sigma^-_{r_p} (H^{(p)})  \leq  \sigma^-_r (H) \! +\! \ep $ and $\sigma^+_{r_p} (H^{(p)})  \geq  \sigma^+_r (H) \! -\! \ep $. Since $\ep$ can be arbitrarily small, we obtain 
$$ \limsup_{p\rightarrow \infty} \sigma^-_r (H^{(p)}) \leq \sigma^-_r (H) \quad \textrm{and} \quad \liminf_{p\rightarrow \infty} \sigma^+_r (H^{(p)}) \geq \sigma^+_r (H) \; , $$
which completes the proof of (\ref{cvtemps}) thanks to (\ref{facidhf}). 

Let us show \eqref{eq: cvTheta}. First, observe that 
\[
\zeta(\Theta_{r_p}(H^{(p)}))=\sigma^+_{r_p}(H^{(p)})-\sigma^-_{r_p}(H^{(p)})\longrightarrow \sigma^+_r(H)-\sigma^-_r(H)=\zeta(\Theta_r(H))
\]
by \eqref{cvtemps}. 
For all $\eta \! \in \! (0, \infty)$, we set $\omega (H, \eta)\! = \! \sup \{ |H_s\! -\! H_{s^\prime}| ; s, s^\prime \! \in \! [0, a ] : |s\! -\! s^\prime| \! \leq \! \eta \}$, the
$\eta$-modulus of uniform continuity of $H$ on $[0, a]$. We have $\lim_{\eta \rightarrow 0} \omega (H, \eta)\! = \! 0$. Since for all $c\in \bbR_+$, $y \! \mapsto y\wedge c$ is $1$-Lipschitz, we get for all $s\! \in \! [0, a]$, 
$$\big| \big( \sigma^-_{r_p} \! (H^{(p)}) \!+\! s\big)\wedge \sigma^+_{r_p} \! (H^{(p)}) - \big( \sigma^-_{r} \! (H) \!+\! s\big)\wedge \sigma^+_{r} \! (H)\big| \! \leq \!  \big| \sigma^-_{r_p} \! (H^{(p)}) \! -\! \sigma^-_{r} \! (H) \big| +\big| \sigma^+_{r_p} \! (H^{(p)}) \! -\! \sigma^+_{r}\!  (H) \big| .$$
Let us denote by $\eta_p$ the number on the right-hand side in the display above. 
By (\ref{cvtemps}), we have $\lim_{p \rightarrow \infty} \eta_p \! = \! 0$. 
Next, we set $y_p:=|\Gamma(H^{(p)})-\Gamma(H)|+|r_p-r|$ and we 
observe that for all $s\! \in \! [0, a]$, 
\[
 \big| \Theta_{r_p} (H^{(p)})(s) \! -\! \Theta_{r} (H)(s) \big|  \leq  \big| H^{(p)}_{ (\sigma^-_{r_p} \! (H^{(p)}) + s)\wedge \sigma^+_{r_p} \! (H^{(p)})} - H_{ (\sigma^-_{r_p} \! (H^{(p)}) + s)\wedge \sigma^+_{r_p} \! (H^{(p)})}\big| + \omega (H, \eta_p) + y_p  \; .
\] 
Thus, 
\[ 
\max_{s\in [0, a]} \big| \Theta_{r_p} (H^{(p)})(s) \! -\! \Theta_{r} (H)(s)\big| \leq \max_{s\in [0, a]} \big| H^{(p)}(s) \! -\! H(s)\big| + \omega (H, \eta_p)+y_p \xrightarrow{p\to\infty} 0\; , 
\]
which completes the proof of the lemma. \cqfd 

\bigskip

\label{pageC}
Let us recall from Section \ref{sec: Levy} the  family of conditional laws $\bN(\, \cdot\, | \, \Gamma(H) \! = \! r)$, $r \! \in \! (0, \infty)$. As implied by Proposition 1.1 of Abraham \& Delmas \cite{AbDe09}, the $\Psi$-height process $H$ under $\bN$ enjoys the following property:  
\[
\forall \, 0< r<u <\infty\,, \quad \Theta_r (H) \; \textrm{ under }\;  \bN\big(\, \cdot\,\big| \, \Gamma(H) =  u\big) \; \eqd \; H \;   \textrm{ under } \; \bN\big(\, \cdot \, \big| \, \Gamma(H)  =  r\big) \; .
\]
Thus, for all nonnegative measurable functional $F\! : \! \bC (\bbR_+, \bbR_+)\! \rightarrow \! \bbR_+$, 
\begin{align*}
\bN \Big[ F \big( \Theta_r (H)  \big)\, \un_{\{ \Gamma(H) >  r \} } \Big] & =  \int_r^\infty \bN \big(\,\Gamma(H)  \in  du\big)  \,\bN \Big[ F \big( \Theta_r (H)  \big) \,\big| \, \Gamma(H)  =   u \Big]  \\
& =   \bN \big(\,\Gamma(H)  >  r\big) \,
\bN \big[ F (H) \, | \, \Gamma(H)  =   r \big] \, , 
\end{align*}
which entails that 
\begin{equation}
\label{tiplawc}
\forall \,  r \in  (0, \infty), \quad \Theta_r (H)  \; \textrm{ under }\; \bN (\, \cdot \, | \, \Gamma(H)  >  r) \; \eqd \; H   \; \textrm{ under } \; \bN (\, \cdot \, | \, \Gamma(H)  =  r) \; .
\end{equation}
Recall that $\bN$-a.e.~for all $s\! \in \! (0, \zeta(H))$, $H_s \! >\! 0$. This also holds true under $\bN (\, \cdot \, | \, \Gamma(H) \! = \! r)$.  Combined with (\ref{tiplawc}), this then implies that $\bN (\, \cdot \, | \, \Gamma(H) \! = \! r)$-a.s.~$H$ satisfies Assumption $(ii)$ of Lemma \ref{conThefun}. 
We also recall that $\bN$-a.e.~there exists a unique time $s\in (0, \zeta(H))$ such that $H_s\! = \! \Gamma(H)$. We readily see that this property still holds true under 
$\bN (\, \cdot \, | \, \Gamma(H) \! = \! r)$. This shows that $\bN(\, \cdot \, | \, \Gamma(H) \! = \! r)$-a.s.~$H$ satisfies Assumption $(i)$ of  Lemma \ref{conThefun}. 

\bigskip

\noi\textbf{Proof of Proposition \ref{bdkjfvs}: }
for all $p\! \geq \! 1$, let $\tau_p\! :\!  \Omega \! \rightarrow \! \bbT_{\! {\rm or}}$ be a GW($\mu_p$)-tree that satisfies the assumptions of Proposition \ref{bdkjfvs}. 
Recall from (\ref{assH}) the sequence $(b_p)_{p\geq1}$ and recall the contour function $(C_s (\tau_p))_{s\in \bbR_+}$ of the ordered rooted tree 
$\tau_p$. We fix $r\! \in (0, \infty) $. To simplify notation, we set $r_p\! = \! \lfloor pr \rfloor/p $. 
The proof of (\ref{heigthcjkb}) given in \cite{Duquesne02} actually 
shows a stronger result: note that the lifetime of $(C_{2b_ps}(\tau_p))_{s\in \bbR_+}$ is equal to $(|\tau_p|\! -\! 1)/b_p$; then the following joint convergence holds weakly on $\bC(\bbR_+, \bbR_+) \! \times \! \bbR_+$ as $p\to \infty$: 
\begin{equation}
\label{jtcontemps}
\Big(\!\big( p^{-1} C_{2b_ps} (\tau_p)\big)_{s\in \bbR+} \, , \, \frac{_{|\tau_p| - 1}}{^{b_p}} \Big)   \textrm{ under }\bP \big(\, \cdot \, \big| \, \Gamma (\tau_p) \! \geq \! pr_p\big)   \xrightarrow{} \big( H , \zeta(H) \big)  \textrm{ under } \bN(\, \cdot \, | \, \Gamma(H) \! \geq  \! r).  
\end{equation}
Indeed, see the proof of Proposition 2.5.2 in \cite{Duquesne02}, page 66. 

By Skorohod's Representation theorem, there exists a probability space $(\Omega^\prime, \cF^\prime, \bP^\prime)$ and processes $H^{(p)}, H \! : \! \Omega^\prime \! \rightarrow \! \bC(\bbR_+, \bbR_+)$ such that: 
\begin{itemize}
\item[$\bullet$] $H^{(p)}$ under $\bP^\prime$ has the same law as $ ( p^{-1}C_{2b_ps} (\tau_p))_{s\in \bbR+}$ under $\bP(\,\cdot\,|\, \Gamma(\tau_p)\ge pr_p)$,  
\vspace{-2mm}
\item[$\bullet$] the law of $H$ under  $\bP^\prime$ is  $\bN (\,\cdot \, | \, \Gamma(H) \! \geq \! r )$,  
\vspace{-2mm}
\item[$\bullet$] $\bP^\prime$-a.s.~$\lim_{p \rightarrow \infty} H^{(p)}\! = \! H$ in $\bC(\bbR_+, \bbR_+)$ and $\lim_{p\rightarrow \infty} \zeta(H^{(p)}) =  \zeta(H)$. 
\end{itemize}
Therefore, $\bP^\prime$-a.s.~$H$ and $H^{(p)}$ satisfy the assumptions of Lemma \ref{conThefun}. Applying the lemma, we get that $\bP^\prime$-a.s.~$\lim_{p\rightarrow \infty} \Theta_{r_p} (H^{(p)}) \! = \! \Theta_r (H)$  in $\bC(\bbR_+, \bbR_+)$. Note that (\ref{tiplawc}) tells that the law of $\Theta_r (H)$ under $\bP^\prime$ is $\bN (\, \cdot \, | \, \Gamma(H) \! =\! r )$. On the other hand, Lemma \ref{tiplaw} and (\ref{contoucsq}) entail that $\Theta_{r_p} (H^{(p)})$ under $\bP^\prime$ has the same law as 
$ (p^{-1}C_{2b_ps} (\tau_p))_{s\in \bbR+}$ under $\bP (\, \cdot \, | \, \Gamma (\tau_p) \!=\! \lfloor pr \rfloor  )$. Indeed, we have shown Proposition \ref{bdkjfvs}. \cqfd 

\section{Proof of Theorem \ref{main}}
\label{mainthpf}

\paragraph{Preliminary computations on GW-trees.} Recall $\Tor$, the set of finite ordered rooted trees.  
Let $t  \in  \Tor$ and $v \in  t$. Recall that $v$ is a tip of $t$ if $|v| =  \Gamma (t)$; that is, $v$ attains the maximum height of $t$. Recall that $k_u(t)$ stands for the number of children of $u$.  
In particular, the children of $\varnothing$ are the single-symbol words $(1), \ldots, (k_\varnothing (t))$.  We also recall $\theta_{(j)} (t)$ from (\ref{thetastem}), the subtree stemming from $(j)$. We introduce the following quantity: 
\begin{equation}
\label{defnunu}
\nu(t)= \big| \big\{ j\in \{ 1, \ldots , k_\varnothing (t) \} : \Gamma \big(\theta_{(j)} (t)\big) \! = \! \Gamma (t)\! -\! 1  \big\} \big|   \; .
\end{equation}
In other words, $\nu(t)$ is the number of individuals in the first generation who have a  tip of $t$ among its descendants.  
Note that $\nu(t) =  0$ if $t$ is reduced to the single vertex $\{ \varnothing\}$; otherwise, 
we always have $1\! \leq \! \nu(t) \! \leq \! k_\varnothing (t) $. 

Let $\mu \! = \! (\mu(k))_{k\in \N}$ be a probability law on $\N$. 
Recall that 
\[ 
g^\mu (s)  \! =\!  \sum_{k\in \mathbb N} s^k \mu (k), \quad s\! \in \! [0, 1] , 
\] 
is the generating function of $\mu$ and that $g^\mu_n$ stands for the $n$-th iteration of $g^\mu$.   
\begin{lem}
\label{nudeux} Let $\mu$ satisfy (\ref{subcrGW}). Let $\tau$ be a GW$(\mu)$-tree on the probability space $(\Omega, \mathcal F, \bP)$. 
Then for all $n\ge 1$, we have
\begin{align}
\label{nudedeux}
\bP \big( \nu(\tau)  \geq  2 \, | \, \Gamma(\tau) =  n\big) &= 1- \frac{g^\mu_{n} (0) \! -\! g^\mu_{n-1} (0)  }{g^\mu_{n+1} (0) \! -\! g^\mu_n (0) } \cdot (g^\mu)^\prime \big( g^\mu_{n-1} (0) \big) \\
\label{nudedeux'}
&\le \frac{1}{(g^\mu)^\prime\big(g^\mu_{n-1}(0)\big)}\bigg(\frac{1-g^\mu_{n-1}(0)}{1-g^\mu_n(0)}+\frac{1-g^\mu_{n+1}(0)}{1-g^\mu_n(0)}-2\bigg).
\end{align}
Here, $(g^\mu)^\prime$ stands for the derivative of $g^\mu$. 
\end{lem}
\noi
\textbf{Proof:} fix $n \geq  1$ and $k \geq  j \geq  2$. Note that 
\[
\bP\big( k_\varnothing(\tau)  =  k; \nu(\tau)  =  j ; \Gamma(\tau)  =  n \big)= \mu (k) \binom{k}{j} \bP\big(\Gamma(\tau) = n-1\big)^j \,\bP\big(\Gamma(\tau) <  n-1\big)^{k-j} \; .
\]
To simplify, we set $b\! := \!\bP (\Gamma (\tau)\!< \! n-1)\! = \! g^\mu_{n-1} (0)$ and $a\! := \! \bP (\Gamma(\tau) \! = \! n-1)\! = \! g^\mu_{n} (0) \! - \! g^\mu_{n-1} (0)$. Then,  
\[
\bP \big(\nu(\tau)  \geq  2 ; \Gamma(\tau)  =  n \big)= \sum_{k\geq 2} \mu (k) \big( (a+b)^k \! -\! b^k \! -\! kab^{k-1} \big) = g^{\mu} (a+b) \! - \! g^\mu (b) \! -\! a (g^\mu)^\prime (b) \; , 
\]
which entails  (\ref{nudedeux}) since $a+b\! = \! g^{\mu}_{n}(0)$. 
As $g^\mu$ is strictly convex under the condition \eqref{subcrGW}, we have $\bP(\Gamma(\tau) \!= \! n)=g^\mu(a+b)-g^\mu(b)\ge a (g^\mu)^\prime(b)$. It follows that  
\begin{align}\notag
\bP \big(\nu(\tau)  \geq   2 \,|\, \Gamma(\tau)  =  n \big)&\le \frac{1}{a(g^\mu)^\prime(b)}\int_b^{a+b} \big((g^\mu)^\prime(u)-(g^\mu)^\prime(b)\big)du \\ \label{nuleqde}
&\le  \frac{(g^\mu)^\prime(a+b)-(g^\mu)^\prime(b)}{(g^\mu)^\prime(b)}=\frac{1}{(g^\mu)^\prime(b)}\int_b^{a+b} (g^\mu)^{\prime\prime}(u)du.
\end{align}
Since $\mu$ is subcritical, $(g^\mu)^\prime (1) \! \leq \! 1$. Combining it with the fact that $(g^\mu)^\prime$ is convex, we get 
\[
\forall \,u \in \big[ b,\, a+b \big], \quad  (g^\mu)^{\prime \prime} (u) \leq \frac{(g^\mu)^{\prime } (1)-(g^\mu)^{\prime } (u)}{1-u} \leq \frac{1-(g^\mu)^{\prime } (u)}{1-(a+b)} \; .
\]
This inequality combined with (\ref{nuleqde}) entails that
\begin{align*}
\bP \big(\nu(\tau)  \geq   2 \,|\, \Gamma(\tau)  =  n \big)  &\leq\frac{1}{(g^\mu)^\prime \big(g^\mu_{n-1}(0)\big)}\cdot \frac{g^{\mu}_{n} (0)-g^{\mu}_{n-1} (0) -g^{\mu}_{n+1} (0) + g^{\mu}_{n} (0)}{1-g^{\mu}_{n} (0)}\\
&= \frac{1}{(g^\mu)^\prime \big(g^\mu_{n-1}(0)\big)}\bigg(\frac{1\! -\! g^\mu_{n-1} (0)}{1\! -\! g^\mu_{n} (0)} + \frac{1\! -\! g^\mu_{n+1} (0)}{1\! -\! g^\mu_{n} (0)} -2\bigg), 
\end{align*}
which is \eqref{nudedeux'}. \cqfd

\begin{figure}[tp]
\centering
\includegraphics[scale=.8]{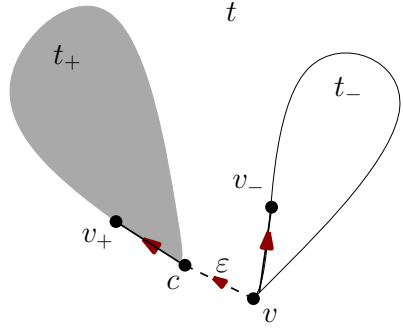}
\caption{\label{fig: rooting} The decomposition of $t$ at the central edge $\varepsilon$ gives rise to two  embedded plane trees $t_-$ and $t_+$. We then root them respectively at the edges $(v, v_-)$ and $(c, v_+)$. Note that during the exploration of $(t, \ep)$ which has been used to define the contour function of $(t, \ep)$, the particle first completes its exploration of  the subtree $t_+$ before it proceeds to $t_-$.  Also, $v_+$ (resp.~$v_-$) is the first edge of $t_+$ (resp.~of $t_-$) visited by the particle. }
\end{figure}

\paragraph{Plane trees viewed from their center, central symmetries.} We discuss here a decomposition of plane trees at their center. Let $t \! = \! (V, E)$ be an embedded plane tree satisfying $|t|>1$. We recall from (\ref{centerdef}) the definition of the center(s) of $t$ and we recall from p.~\pageref{p: ed}  
the definition of the central edges of $t$. Note that central edges are oriented edges. 
Let $\ep \! =\!  (v, c) \! \in \! K(t)$. Observe that $c$ is necessarily a center of $t$. The removal of $\ep$ splits $t$ into two subtrees $t_-$ and $t_+$: $t_-$ being the one that contains $v$ and $t_+$ the one that contains $c$. Both are embedded plane trees. We root them in the following way (see also Figure \ref{fig: rooting}). 
\vspace{1mm}
\begin{compactenum}[-]
\item[$\bullet$] Let $v_-$ be a neighbor of $v$ such that $v_-$ is next to $c$ in the cyclic order on the set of neighbors of $v$ induced by the orientation of the plane. Note that $(v, v_-)$ is an oriented edge of $t_-$. As explained in Section \ref{sec: pl}, 
the edge-rooted plane tree $(t_-, (v,v_-))$ induces an ordered rooted tree $(t_-, (v,v_-))_o$ that we denote by $T_- (t, \ep)   \in  \Tor$ in what follows.
\item[$\bullet$] Let $  v_+$ be the neighbor that is next to $v$ in the cyclic order on the set of the neighbors of $c$ induced by the orientation of the plane. Note that $(c, v_+)$ is an oriented edge of $t_+$. 
Then, $(t_+, (c,v_+))$ induces an ordered rooted tree $(t_+, (c,v_+))_o$ that we denote by $T_+ (t, \ep)   \in  \Tor$ in what follows. 
\end{compactenum}
\vspace{0.5mm}
It is important to note that if $(t, \ep)$ and $(t', \ep')$ are two equivalent edge-rooted plane trees, then $T_- (t, \ep) \! =\! T_- (t^\prime, \ep^\prime) $ and $T_+ (t, \ep) \! =\! T_+ (t^\prime, \ep^\prime)$. This shows that $T_-(t, \ep)$ and $T_+(t, \ep)$ only depend on the equivalence class of the edge-rooted plane tree $(t, \ep)$. Then, they only depend on the ordered rooted tree $(t, \ep)_o$. For this reason, we sometimes write $T_{+/-}((t, \ep)_o)$ instead of $T_{+/-}(t, \ep)$. 

Recall from (\ref{defnunu}) the definition of $\nu$. Let us observe the following. 
\begin{align}
\label{descri}
&\textrm{If $\diam (t) =  2p\!+\!1$, $p\!\in\! \mathbb N$, then $\Gamma (T_-(t, \ep)) =  \Gamma (T_+ (t, \ep)) =  p$ and $|K(t)|=  2$.}\\
\label{descrii}
&\textrm{If $\diam (t) =  2p$, $p\!\in\! \mathbb N$, then $\Gamma (T_-(t, \ep)) =  p\! -\! 1$, $\Gamma (T_+ (t, \ep)) =  p$ and $|K(t)| = 1\!+\! 
\nu (T_+(t, \ep))$.}
\end{align}
We introduce the \textit{number of central symmetries} of $(t, \ep)$ as follows: 
\begin{equation}
\label{def: Sym}
\mathtt{Sym} (t,\ep)= \big| \big\{ \ep^\prime \! \in \! K(t) : (t, \ep^\prime)_o \! = \! (t, \ep)_o\big\} \big| \; .
\end{equation}
Note that $ \mathtt{Sym} (t,\ep)$ only depends on the equivalence class of the edge-rooted tree $(t, \ep)$. So we may sometimes write 
$ \mathtt{Sym} ((t,\ep)_o)$ instead of $ \mathtt{Sym} (t,\ep)$.

Let $p$ be a positive integer. Recall from (\ref{def: Tplk}) the set $\bbT_{\! {\rm pl}}(p)$ of plane trees with diameter $p$. 
We introduce the following notation. 
\begin{equation}\label{def: B_p}
B_p \! := \! \big\{ (t,\ep)\, : \;  t \! \in \!  \Tpl(p),  \; \ep \! \in \! K(t)  \big\} \quad \textrm{and} \quad  B^o_p \! := \! \big\{ (t,\ep)_o \! \in \! \Tor \, :\;  (t, \ep) \! \in \! B_p   \big\} \; .
\end{equation}
If we denote $\pi: (t, \ep)\in B_p \mapsto (t, \ep)_o$ the canonical projection, then $\pi^{-1}((t, \ep)_o)$ contains a number $\mathtt{Sym}((t,\ep)_o)$ of elements, for each $(t, \ep)_o\in B_p^o$. 
Recall that $\Tor^{^=} (p) =  \{ t\in \Tor: \Gamma (t)  =  p\}$. 
It is not difficult to check that 
\begin{equation}
\label{concombi}
\textrm{the mapping} \;  (t, \ep)_o \! \in \!  B^o_p \longmapsto \big( T_-(t, \ep), T_+ (t, \ep)  \big) \in \Tor^{^=} (\lfloor \tfrac{p-1}{2} \rfloor  )\times  \Tor^{^=} (\lfloor \tfrac{p}{2} \rfloor ) \; \textrm{is bijective.}
\end{equation}
As already mentioned, $ \mathtt{Sym} (t,\ep)$ only depends on $(t, \ep)_o$. Combined with  (\ref{concombi}), we see that $ \mathtt{Sym} (t,\ep)$ is a function of  $\big(T_-(t, \ep), T_+ (t, \ep)  \big)$. Denote by $S_p$ this function: $S_p$ is the unique function 
from $\Tor^{^=} (\lfloor \tfrac{p-1}{2}\rfloor )\times  \Tor^{^=} (\lfloor \frac{p}{2} \rfloor )$ to $\{ 1, 2, \ldots\}$ such that 
\begin{equation}
\label{esspedef}
\forall  (t, \ep) \in B_p, \quad  \mathtt{Sym} (t,\ep)= S_p \big( T_-(t, \ep), T_+ (t, \ep)  \big) \; .
\end{equation}
Let us briefly discuss the properties of $S_p$. We first consider the case of plane trees with an odd diameter. 
Let $t\! \in \! \Tpl (2p\!+\!1)$. The tree $t$ has two central edges $\ep \! := \!  (c, c^\prime)$ and $\ep^\prime \! := \! (c^\prime, c)$. 
Then $\mathtt{Sym} (t, \ep)\! = \! 2$ if and only if $(t, \ep)$ and $(t, \ep^\prime)$ are equivalent. In this case, $T_-(t, \ep)\! = \! T_+(t, \ep)$. We have the following   
\begin{equation}
\label{odddia} \forall t_1, t_2 \! \in \!  \Tor^{^=}(p), \quad S_{2p+1} (t_1, t_2) = 1+ \un_{\{ t_1= t_2\}} \; .
\end{equation}
We next consider the case of plane trees with an even diameter. Let $t\! \in \! \Tpl(2p)$. Then $t$ has a unique center $c$. Let us write $t_1\! = \! T_-(t, \ep)$ and $t_2\! = \! T_+(t, \ep)$. We also denote $N=k_\varnothing(t_2)$, the number of children of the root of $t_2$. 
Note that $\degree(c) \! = \! 1+ N$. Recall that $\theta_{(j)} (t_2)$ stands for the subtree stemming from the $j$-th child of the root of $t_2$, for $1\le j\le N$. 
Then,  the number of internal symmetries is $\degree(c)/d$, where $d$ is the minimal period of the list $(\theta_{(1)} (t_2), \ldots,  \theta_{(N)} (t_2) , t_1) $. However, we will only need the following bound: 
\begin{equation}
\label{evendia} \forall (t_1, t_2)\! \in \! \Tor^{^=} (p\! -\! 1) \! \times \!  \Tor^{^=} (p), \;    S_{2p} (t_1, t_2) \leq 1+ \big| \big\{ j\! \in \! \{1, \ldots, N  \}: t_1\! = \! \theta_{(j)} (t_2)\big\}   \big|  \!  \leq \!   
1+ \nu (t_2)  .
\end{equation}
The last inequality comes from the fact that $\Gamma (t_1)\! = p - 1 =  \Gamma (t_2) - 1$.

\medskip

Let $\mu \! = \! (\mu(k))_{k\in \N}$ satisfy (\ref{subcrGW}).
For a plane tree $t\in\Tpl$, we recall from (\ref{jpkgkjh}) the weight $W_{\! \mu}(t)$.
We define  
$$ \forall\, t \! \in \! \bbT_{\! {\rm pl}}, \; \forall \ep \! \in \! K(t), \quad  \Wt (t, \ep)\! := \!  W_\mu (t) / |K(t)| \; .$$
Observe that $ \Wt (t, \ep)$ only depends on the induced ordered rooted tree $(t, \ep)_o $. For this reason, we will write $\Wt ((t, \ep)_o)$ instead of $ \Wt (t, \ep)$. Let $\tau$ and $\widetilde{\tau}$ be two independent GW($\mu$)-trees defined on $(\Omega, \cF, \bP)$.
We note that 
\begin{eqnarray}
\label{gerboulade}
 \Wt (t, \ep) & =&  |K(t)|^{-1} \!\!\!\!\!\prod_{\quad v\in T_-(t, \ep)}   \mu \big( \degree (v) \! -\! 1\big) \!   \prod_{\quad v\in T_+(t, \ep)}   \mu \big( \degree (v) \! -\! 1\big) \nonumber \\
 & =&  |K(t)|^{-1} \; \bP \big(\tau=T_-(t, \ep) \,; \widetilde{\tau}=T_+(t, \ep) \big) \,  . 
 \end{eqnarray}
Recall from \eqref{def: B_p} the sets $B_p$ and $B_p^o$. 
By definition, we have $\Wt(t, \ep)=\Wt(t, \ep')$, if both $\ep, \ep'\in K(t)$. 
Also, recall from \eqref{def: Sym} the number $\mathtt{Sym}((t, \ep)_o)$. 
Then,
\[
Z_{p} (\mu)  = \sum_{t\in \Tpl(p)} W_\mu(t)  \; = \sum_{(t, \ep) \in B_{p}} \Wt  (t, \ep) \;=  \sum_{(t, \ep)_ o \in B^o_{p}}  \mathtt{Sym} ((t, \ep)_o) \Wt  ((t, \ep)_o)\,. 
\]
Recall that if $\diam(t)$ is even, then $|K(t)| \! = \! 1+ \nu (T_+(t, \ep))$. Recall from \eqref{esspedef} the function $S_p$. 
It follows from the above display, \eqref{concombi} and \eqref{gerboulade} that 
\begin{equation}
\label{floupii}
Z_{2p}(\mu)= \!\! \!\! \sum_{t_1 \in  \Tor^{^=}(p-1)}\sum_{t_2 \in  \Tor^{^=}(p)} \!\! \!\!\bP \big(\tau=t_1;\widetilde{\tau}=t_2 \big) \frac{S_{2p} (t_1, t_2)}{1+ \nu(t_2)}=  \! \bE \Big[ \frac{S_{2p} (\tau, \widetilde{\tau}) }{1+ \nu (\widetilde{\tau}) }\un_{\{ \Gamma (\tau)+1= \Gamma (\widetilde{\tau}) = p \}}  \Big] \,.
\end{equation}
We obtain in a similar way that   
\begin{equation}
\label{floupi}
Z_{2p+1} (\mu)  =  \bE \big[ \tfrac{1}{2}S_{2p+1} (\tau, \widetilde{\tau}) \un_{\{ \Gamma (\tau)=  \Gamma (\widetilde{\tau}) = p \}}  \big]\; .
\end{equation}
From (\ref{floupii}) and (\ref{evendia}), we deduce that  $Z_{2p} (\mu) \! \le \!  \bP (\Gamma (\tau) \! =\! p\! -\! 1; \Gamma (\widetilde{\tau}) \! =\!  p)$. 
Similarly, (\ref{floupi}) and (\ref{odddia}) implies that $Z_{2p+1} (\mu) \! \le \!  \bP (\Gamma (\tau) \! =\! \Gamma (\widetilde{\tau})  \! =\! p)$. We have proved the following lemma. 
\begin{lem}
\label{Zpfini} Let $\mu$ satisfy (\ref{subcrGW}). Let $Z_p$ be as defined in  (\ref{jpkgkjh}). Then, $Z_p (\mu) <  \infty$, for all $p\in\mathbb N$. 
\end{lem}  

It follows from Lemma \ref{Zpfini} that  for all $t\! \in \! \bbT_{\! {\rm pl}} (p)$, $Q^\mu_p (t)\! = \! W_{\! \mu} (t)/ Z_p (\mu)$ is well defined. We next show the following. 
\begin{lem}
\label{splitter} Let $\mu$ satisfy (\ref{subcrGW}). Let $p\! \geq \! 1$ be an integer and let $(\rT_p, \bep_p)$ be a pair of random variables defined on $(\Omega, \mathcal F, \bP)$ such that
\[
\forall (t, \ep)\! \in \! B_p\, , \quad \bP \big((\rT_p, \bep_p)\! = \! (t, \ep) \big)\! = \! Q^\mu_p (t)/ |K(t)|= Z_p(\mu)^{-1}\Wt (t, \ep) \; .
\]
Let $\tau, \widetilde{\tau}\! : \! \Omega \!  \rightarrow \! \bbT_{\! {\rm or}}$ be two independent GW($\mu$)-trees. 
Let $F, G \! : \! \bbT_{\! {\rm or}} \! \rightarrow \! \bbR_+$ be two bounded nonnegative measurable functions. Then, the following holds true: 
\begin{equation}
\label{zoicoic}
\bE \big[ F\big( T_-(\rT_{2p+1}, \bep_{2p+1}) \big) G\big( T_+(\rT_{2p+1}, \bep_{2p+1}) \big) \big] = \frac{ \bE \big[ S_{2p+1}( \tau, \widetilde{\tau}) F(\tau) G(\widetilde{\tau}) \un_{\{ \Gamma (\tau) = \Gamma (\widetilde{\tau})= p\}} \big] }{\bE \big[ S_{2p+1}( \tau, \widetilde{\tau})  \un_{\{ \Gamma (\tau) = \Gamma (\widetilde{\tau})= p\}} \big]}. 
\end{equation}
\vspace{-1mm}
\begin{equation}
\label{zoicouic}
\bE \big[ F\big( T_-(\rT_{2p}, \bep_{2p}) \big) G\big( T_+(\rT_{2p}, \bep_{2p}) \big) \big] = \frac{\bE \Big[ \frac{{S_{2p} (\tau, \widetilde{\tau}) }}{{1+ \nu (\widetilde{\tau}) }} F(\tau) G(\widetilde{\tau})  \un_{\{ \Gamma (\tau)+1 = \Gamma (\widetilde{\tau}) = p \}}  \Big]}{\bE \Big[ \frac{{S_{2p} (\tau, \widetilde{\tau}) }}{{1+ \nu (\widetilde{\tau}) }}\un_{\{ \Gamma (\tau)+1= \Gamma (\widetilde{\tau}) = p \}}  \Big]}\; .
\end{equation}
\end{lem}
\noi
\textbf{Proof:} by definition, we have
\begin{align*}
&\,\bE  \Big[  F\big( T_- (\rT_{2p}, \bep_{2p}) \big) G\big( T_+ (\rT_{2p}, \bep_{2p})  \big)  \Big] \\ 
= &\; \tfrac{1}{Z_{2p}(\mu)}  \sum_{(t, \ep) \in B_{2p}}  \Wt  (t, \ep) F\big(T_-(t, \ep)\big) G\big(T_+(t, \ep)\big)  \\
 = &\; \tfrac{1}{Z_{2p}(\mu)}  \sum_{(t, \ep)_o \in B^o_{2p}} \mathtt{Sym}  \big( (t, \ep)_o) \Wt   ( (t, \ep)_o \big) F\big(T_-((t, \ep)_o) \big) G\big(T_+((t, \ep)_o) \big),
\end{align*}
where we recall from \eqref{def: Sym} the definition of $ \mathtt{Sym}$. Applying \eqref{concombi} and then the same argument in \eqref{floupii}, we find that 
\begin{align*}
&\,\bE  \Big[ F\big( T_- (\rT_{2p}, \bep_{2p}) \big) G\big( T_+ (\rT_{2p}, \bep_{2p})  \big)  \Big]\\
 = &\;\frac{1}{Z_{2p}(\mu)} \sum_{t_1 \in  \Tor^{^=}(p-1)}\sum_{t_2 \in  \Tor^{^=}(p)}  
\bP(\tau=t_1; \widetilde\tau=t_2)\frac{S_{2p} (t_1, t_2)}{1+ \nu(t_2)} F(t_1) G(t_2)  \\
= &\;\frac{1}{Z_{2p}(\mu)} \, \bE \bigg[ \frac{S_{2p} (\tau, \widetilde{\tau}) }{1+ \nu (\widetilde{\tau}) } F(\tau) G(\widetilde{\tau})  \un_{\{ \Gamma (\tau)+1 = \Gamma (\widetilde{\tau}) = p \}}  \bigg] \; , 
\end{align*}
which entails (\ref{zoicouic}) by (\ref{floupii}). We prove (\ref{zoicoic}) in a similar way. \cqfd

\begin{lem}
\label{contigu} Let $\mu$ satisfy (\ref{subcrGW}). Let $p\ge 1$ be an integer  and let $(\rT_p, \bep_p)$ be the pair as in Lemma \ref{splitter}. 
Let $\tau, \widetilde{\tau} :  \Omega  \rightarrow \Tor$ be two independent GW($\mu$)-trees. 
We set 
\begin{align*}
a_p (\mu) := \,& \bP \Big( S_p(\tau,  \widetilde{\tau} ) \geq  2\, \Big|\,  \Gamma(\tau_ o) =  \lfloor \tfrac{1}{2} (p - 1) \rfloor ; \Gamma( \widetilde{\tau})  =  \lfloor \tfrac{1}{2} p \rfloor \Big)\,; \\
b_p (\mu)  :=\,&  \bP \Big( \nu (\tau) \! \geq \! 2 \,\Big| \,\Gamma (\tau)   = \lfloor \tfrac{1}{2}  p \rfloor \Big). 
\end{align*}
Let $F, G  :  \Tor \rightarrow \bbR_+$ be two bounded measurable functions. Then, the following holds true: 
\begin{align}
\label{concontig}
&\Big| \bE \Big[ F\big( T_-(\rT_{p}, \bep_{p}) \big)\cdot G\big( T_+(\rT_{p}, \bep_{p}) \big) \Big]- \bE \Big[ F(\tau) \,\Big|\, \Gamma (\tau ) =  \lfloor \frac{_{_1}}{^{^2}} (p - 1) \rfloor\Big] \cdot
\bE \Big[ G(\tau) \,\Big|\, \Gamma (\tau ) =  \lfloor  \frac{_{_1}}{^{^2}} p \rfloor \Big] \Big| \nonumber \\  
&\qquad\qquad\qquad  \leq 4 \lVert F \rVert_{\infty}  \lVert G \rVert_{\infty}  \big(a_p (\mu)+b_p (\mu) \big) \; .
\end{align}
\end{lem}
\noi
\textbf{Proof:} we only detail the case of even diameters. The case of odd diameters can be treated similarly. To ease the writing, we set 
\begin{align*}
&\qquad \qquad\qquad \qquad  \alpha_{2p}(F, G) =  \bE \Big[\frac{_{S_{2p} (\tau, \widetilde{\tau}) -1}}{^{1+ \nu (\widetilde{\tau}) }} F(\tau) G(\widetilde{\tau})  \un_{\{ \Gamma (\tau)+1 = \Gamma (\widetilde{\tau}) = p \}}  \Big], \; \\
\beta_{2p}& = \bE \Big[ \frac{_{1}}{^{1+ \nu (\widetilde{\tau}) }} F(\tau) G(\widetilde{\tau})  \un_{\{ \Gamma (\tau)+1 = \Gamma (\widetilde{\tau}) = p \}}  \Big],  \quad
\gamma_{2p}(F, G)=  \tfrac{1}{2}\,\bE \Big[ F(\tau) G(\widetilde{\tau})  \un_{\{ \Gamma (\tau)+1 = \Gamma (\widetilde{\tau}) = p \}}  \Big].
\end{align*}
First, we note that 
\begin{eqnarray*}
 |\alpha_{2p}(F, G)\! -\! \beta_{2p}| & \leq & \bE \Big[ \frac{_{S_{2p} (\tau, \widetilde{\tau}) -1}}{^{1+ \nu (\widetilde{\tau}) }}\,\un_{\{ S_{2p} (\tau, \widetilde{\tau}) \geq 2 \}} \big| F(\tau) G(\widetilde{\tau}) \big|\,  \un_{\{ \Gamma (\tau)+1 = \Gamma (\widetilde{\tau}) = p \}}  \Big] \\
 & \leq & \lVert F \rVert_{ \infty}   \lVert G \rVert_{\infty}\,  a_{2p}(\mu) \, \bP \big( \Gamma (\tau)+1  = \Gamma (\widetilde{\tau})  =  p \big),
\end{eqnarray*} 
since $S_{2p} (\tau, \widetilde{\tau})\! \leq \! 1+ \nu (\widetilde{\tau})$ by (\ref{evendia}). Next, we observe that 
\begin{eqnarray*}
 |\beta_{2p} \! -\! \gamma_{2p}(F, G)| & \leq & \bE \Big[ \Big( \tfrac{1}{2}-\frac{_{1}}{^{1+ \nu (\widetilde{\tau}) }} \Big)
 \un_{\{ \nu ( \widetilde{\tau}) \geq 2 \}} \big| F(\tau) G(\widetilde{\tau}) \big| \un_{\{ \Gamma (\tau)+1 = \Gamma (\widetilde{\tau}) = p \}}  \Big] \\
 & \leq & \lVert F \rVert_{ \infty}   \lVert G \rVert_{\infty}\,   b_{2p}(\mu)\, \bP \big( \Gamma (\tau)+1  =  \Gamma (\widetilde{\tau})  =  p \big)\; .
\end{eqnarray*} 
Thus, 
\begin{equation}
\label{bd_alpha}
 |\alpha_{2p}(F, G) - \gamma_{2p}(F, G)| \leq  \lVert F \rVert_{ \infty}   \lVert G \rVert_{\infty}\,  \big( a_{2p}(\mu) +b_{2p}(\mu) \big) \bP \big( \Gamma (\tau)+1  =  \Gamma (\widetilde{\tau})  =  p \big) \; .
 \end{equation}
Write $\mathbf 1$ for the constant function $F\equiv 1$. By \eqref{zoicouic}, we have
\[
\bE \big[ F\big( T_-(\rT_{2p}, \bep_{2p}) \big) G\big( T_+(\rT_{2p}, \bep_{2p}) \big) \big]=\frac{\alpha_{2p}(F, G)}{\alpha_{2p}(\mathbf 1, \mathbf 1)}. 
\]
Then,  
\begin{align*}
&\ \Big| \bE \big[ F\big( T_-(\rT_{2p}, \bep_{2p}) \big) G\big( T_+(\rT_{2p}, \bep_{2p}) \big) \big]- \bE \big[ F(\tau) \,\big|\, \Gamma (\tau )\! =  p-1 \big] \cdot
 \bE \big[ G(\tau) \,\big|\, \Gamma (\tau )\! = p  \big] \Big|  \\ 
 \le &\  \frac{\alpha_{2p}(F, G)}{\alpha_{2p}(\mathbf 1, \mathbf 1)\gamma_{2p}(\mathbf 1, \mathbf 1)}\,\Big| \gamma_{2p}(\mathbf 1, \mathbf 1)-\alpha_{2p}(\mathbf 1, \mathbf 1)\Big|+\frac{1}{\gamma_{2p}(\mathbf 1, \mathbf 1)}\,\Big|\alpha_{2p}(F, G)-\gamma_{2p}(F, G)\Big|,
\end{align*}
which entails (\ref{concontig}) by \eqref{bd_alpha} in the case of an even diameter. \cqfd 

\begin{lem}
\label{symevava}  Let $\mu$ satisfy (\ref{subcrGW}). Let $p\ge 1$ be an integer. 
Let $\tau, \widetilde{\tau}:  \Omega  \rightarrow \Tor$ be two independent GW($\mu$)-trees. 
Let the function $S_p$ be defined as in (\ref{esspedef}). Then, 
\begin{equation}
\label{beurkhhh}
\bP \Big( S_p(\tau,  \widetilde{\tau} )  \geq  2  \,\big|\,  \Gamma(\tau)  =  \lfloor \frac{_{_1}}{^{^2}} (p - 1) \rfloor ; \Gamma( \widetilde{\tau})  =  \lfloor \frac{_{_1}}{^{^2}} p \rfloor \Big)
  \leq \frac{
\bP \big( \tau =  \widetilde{\tau}  \,\big|\,  \Gamma(\tau) =  \Gamma( \widetilde{\tau})  =  \lfloor \tfrac{1}{2} (p - 1) \rfloor  \big)}{(g^\mu)^\prime\big(g^\mu_{p-1}(0)\big)} \; .
\end{equation}
\end{lem}
\noi
\textbf{Proof:} in the odd diameter cases, \eqref{beurkhhh} is a combined consequence of  (\ref{odddia}) and the fact that $0<(g^\mu)^\prime(u)\le 1$ for all $u\in[0, 1]$, as $\mu$ is (sub)critical. Let us consider the even diameter case. We apply (\ref{evendia}) to get the following.
\begin{align*}
\bP \big( S_{2p}(\tau,  \widetilde{\tau} )  \geq  2 \,;  \Gamma(\tau)+1 = \Gamma( \widetilde{\tau}) \! = \! p \big) 
 &\leq   \bE \Big[ \sum_{1\leq j\leq k_\varnothing ( \widetilde{\tau})}  \un_{\{\theta_{(j)}  \widetilde{\tau}\,=\, \tau   \}\,\cap\, \{ \Gamma(\tau)\,=\,p-1\} }\,\un_{\{\Gamma(\widetilde\tau)\,=\,p\}} \Big] \\
& \leq  \sum_{k\in \N } k\mu(k)\, \bP \big(\widetilde{\tau} \! = \! \tau \,;  \Gamma(\tau) \! =\! p-1 \big)\,\bP\big(\,\Gamma(\tau)\le p-1\big)^{k-1} \\
& \leq \bP \big(\widetilde{\tau} \! = \! \tau \,;  \Gamma(\tau)  = p-1 \big)(g^\mu)^\prime\big(g^\mu_{p}(0)\big). 
\end{align*}
Recall that $\bP(\Gamma( \widetilde{\tau}) \! = \! p )=g^\mu_{p+1}(0)-g^\mu_p(0)$. We then deduce that
\begin{align*}
&\bP \big( S_{2p}(\tau,  \widetilde{\tau} )  \geq  2 \, \big|\,  \Gamma(\tau)+1  =  \Gamma( \widetilde{\tau})  =  p \big)\\
&\qquad\qquad\qquad \le \bP \big(\widetilde{\tau}  =  \tau \,\big|\,  \Gamma(\tau)=\Gamma(\widetilde{\tau})  = p-1 \big)(g^\mu)^\prime\big(g^\mu_{p}(0)\big)\cdot \frac{g^\mu_p(0)-g^\mu_{p-1}(0)}{g^\mu_{p+1}(0)-g^\mu_p(0)}\; ,
\end{align*}
 which implies \eqref{beurkhhh}, since $g^\mu_{p+1}(0)-g^\mu_p(0)\ge (g^\mu)^\prime(g^\mu_{p-1}(0))(g^\mu_p(0)-g^\mu_{p-1}(0))$ and $(g^\mu)^\prime(g^\mu_{p}(0))\le 1$, under the assumption \eqref{subcrGW}.
\cqfd

\bigskip

\paragraph{Main proof.} Recall that $\bC (\bbR_+, \bbR_+)$ is equipped with the Polish topology of the uniform convergence on the compact subsets. 
Recall $\bC$ from p.~\pageref{pageC}, the set of coding functions. In particular, if $H\! \in \! \bC$, its lifetime $\zeta(H) \! = \! \sup \{ t\in \bbR_+ : H_t \! >\! 0 \}\in (0, \infty)$. 
 Recall from (\ref{concadef}) the concatenation $H\oplus \widetilde{H}$ of two coding functions $H$ and $\widetilde H$.  
We need the following lemma whose proof is direct (and is thus omitted). 
\begin{lem}
\label{concatef} Let $(H^{(p)}, p\in \N)$ and $(\widetilde{H}^{(p)}, p\in \N)$ be two sequences of coding functions. Let $H, \widetilde H\in \bC$. 
W assume the following conditions.  
\vspace{-2mm}
\begin{itemize}
\item[(i)] $ \lim_{p \rightarrow \infty} H^{(p)}\! = \! H $ and $\lim_{p \rightarrow \infty} \widetilde{H}^{(p)}\! = \! \widetilde{H}\, $ in $\bC(\bbR_+, \bbR_+)$. 
\vspace{-2mm}
\item[(ii)]  $\lim_{p\rightarrow \infty} \zeta(H^{(p)})\! = \! \zeta(H)$ and  $ \lim_{p \rightarrow \infty} \zeta(\widetilde{H}^{(p)}) \! = \! \zeta(\widetilde{H})$. 
\end{itemize}
\vspace{-2mm}
Then,  
\begin{equation}
\label{cheigconca}
\lim_{p \rightarrow \infty} \Gamma (H^{(p)}) \! = \! \Gamma (H) \quad \textrm{and} \quad \lim_{p\rightarrow \infty } H^{(p)} \! \oplus \!  \widetilde{H}^{(p)}\! =\! H\oplus \!  \widetilde{H} \; \textrm{ in $\bC(\bbR_+, \bbR_+)$.} 
\end{equation}
\end{lem}

We apply (\ref{cheigconca}), Lemma \ref{nudeux} and the convergence in (\ref{jtcontemps}) to show the following. 
\begin{lem}
\label{nupzer} Let $(\mu_p, p\! \geq \! 1)$ be a sequence of laws that satisfy \eqref{subcrGW}. 
Assume that (\ref{assH}) and (\ref{asszeta}) take place. For all $p\! \geq \! 1$, let $\tau_p  : \Omega  \rightarrow  \Tor$ be a GW($\mu_p$)-tree. Let us fix $r\in  (0, \infty)$. 
Then, 
\[ 
\bP \big( \, \nu( \tau_p)  \geq  2  \,  \big| \, \Gamma (\tau_p) =  \lfloor pr \rfloor  \big) \xrightarrow{p\to\infty} 0 \; .
\]
\end{lem}
\noi
\textbf{Proof:} write $r_p:=\lfloor pr \rfloor$. By (\ref{jtcontemps}) and the first limit in  (\ref{cheigconca}), 
we get 
\[
\tfrac{1}{p}\, \Gamma (\tau_p)   \  \textrm{ under }\  \bP (\, \cdot \, | \, \Gamma (\tau_p)  \geq  pr_p  ) \  \xrightarrow[p\to\infty]{\text{(d)}} \ \Gamma(H) \ \textrm{ under }\ \bN(\, \cdot \, | \, \Gamma(H)  \geq   r)\,.  
\]
Since the law of $\Gamma(H)$ under $\bN(\, \cdot \, | \, \Gamma(H) \! \geq  \! r)$ is diffuse, for all $s\! \in \! [r, \infty)$ and for all nonnegative integers $p_0$, we deduce the following convergence.
\begin{equation}
\label{flkbe}
\frac{1-g^{\mu_p}_{\lfloor ps\rfloor+p_0}(0)}{1-g^{\mu_p}_{\lfloor pr\rfloor}(0)}=
\bP \big( \,\Gamma (\tau_p) \geq  \lfloor ps \rfloor + p_0 \,  \big| \, \Gamma (\tau_p) \geq  \lfloor pr \rfloor  \big) \; \xrightarrow[p\rightarrow \infty] \;  \bN ( \,\Gamma(H)   \geq  s \, | \, \Gamma(H) \geq  r) \; , 
\end{equation}
where we have used the fact that $\bP ( \Gamma (\tau_p) \ge n) = 1 - g^{\mu_p}_n (0)$, 
recalling that $g^{\mu_p}_n$ stands for the $n$-th iteration of the generating function $g^{\mu_p}$ of $\mu_p$. Then, (\ref{flkbe}) entails that 
\begin{equation}
\label{bd1}
 \forall \, s  \in  (0, \infty), \; \forall \, p_0  \in  \N, \quad \frac{1-g^{\mu_p}_{\lfloor ps \rfloor + p_0} (0)}{1-g^{\mu_p}_{\lfloor ps \rfloor } (0)}\; \xrightarrow[p\rightarrow \infty] \; 1 \; .
\end{equation} 
On the other hand, note that for $u\in (0, 1)$, $\lfloor pu\rfloor \le p-1$ for sufficiently large $p$. Then by the convexity of $g^{\mu_p}$, we have
\begin{equation}
\label{bd2}
\big(g^{\mu_p}\big)^\prime \big(g^{\mu_p}_{p-1}(0)\big)\ge \frac{g^{\mu_p}\big(g^{\mu_p}_{p-1}(0)\big)-g^{\mu_p}\big(g^{\mu_p}_{\lfloor pu \rfloor}(0)\big)}{g^{\mu_p}_{p-1}(0)-g^{\mu_p}_{\lfloor pu \rfloor}(0)}=\frac{\frac{1-g^{\mu_p}_{\lfloor pu\rfloor+1}(0)}{1-g^{\mu_p}_{\lfloor pu\rfloor}(0)}-\frac{1-g^{\mu_p}_{p}(0)}{1-g^{\mu_p}_{\lfloor pu\rfloor}(0)}}{1-\frac{1-g^{\mu_p}_{p-1}(0)}{1-g^{\mu_p}_{\lfloor pu\rfloor}(0)}}\xrightarrow{p\to\infty} 1,
\end{equation}
by \eqref{flkbe}. 
We deduce from this, \eqref{bd1} and Lemma \ref{nudeux} the desired result. \cqfd 

\begin{lem}
\label{zouloufli} Let $(\mu_p, p\! \geq \! 1)$ be a sequence of laws that satisfy \eqref{subcrGW}. 
Assume that (\ref{assH}) and (\ref{asszeta}) take place. For all $p\! \geq \! 1$, let $\tau_p  : \Omega  \rightarrow  \Tor$ be a GW($\mu_p$)-tree. Let us fix $r\in  (0, \infty)$ and set $r_p=\lfloor pr \rfloor$.  
Then, 
\begin{align}
\label{nul_ev1}
&\bP \big( \tau_p=  \widetilde{\tau}_p  \,\big|\,  \Gamma(\tau_p)  =  \Gamma( \widetilde{\tau}_p)   =  r_p  \big)\xrightarrow{p\to\infty} 0 \; .\\
\label{nul_ev2}
&\bP \big( S_p(\tau_p,  \widetilde{\tau}_p )  \geq  2 \, \big|\,  \Gamma(\tau_p) =  \lfloor \frac{_{_1}}{^{^2}} (r_p - 1) \rfloor ; \Gamma( \widetilde{\tau})  =  \lfloor \frac{_{_1}}{^{^2}} r_p \rfloor \big) \xrightarrow{p\to\infty} 0 \; . 
\end{align}
\end{lem}
\noi
\textbf{Proof:} let $H, \widetilde{H}\! :\!  \Omega \! \rightarrow \bC (\bbR_+, \bbR_+)$ be two independent processes with common law $\bN (\, \cdot \, \big| \, \Gamma(H) \! = \! r)$.
By (\ref{heightcon}) in Proposition \ref{bdkjfvs}, we obtain the following weak convergence on $\bbR_+^2$: 
\begin{equation}
\label{dflkbv}
\big( \frac{{_1}}{^{b_p}} | \tau_p| \, , \,\frac{{_1}}{{^{b_p}}}  |\widetilde{\tau}_p |  \big) \quad \textrm{under $\bP (\, \cdot \, |\,  \Gamma(\tau_p)  =  \Gamma( \widetilde{\tau}_p)   = \lfloor pr \rfloor  )$} \; \overset{{\rm (d)}}{\underset{p\rightarrow \infty}{-\!\!\! -\!\!\! -\!\!\! \longrightarrow}} \;  \big( \zeta (H) , \zeta( \widetilde{H}) \big) \; .
\end{equation}
Let $\Delta \! = \! \{ (x,x) ;  x\! \in \! \bbR_+ \}$ be the diagonal of $\bbR_+^2$. The distribution of $\zeta(H)$ under $\bN (\, \cdot\, \big| \, \Gamma(H) \! = \! r)$ is diffuse. It follows that  
$\bP (  ( \zeta (H) , \zeta( \widetilde{H}) )\in \Delta )\! = \!  \bP ( \zeta (H) \! = \!  \zeta( \widetilde{H}) )\! = \! 0$. Since $\Delta$ is a closed set, applying Portmanteau's Theorem (see for instance Ethier \& Kurtz \cite{EthKur86}, Theorem 3.1 $(a) \! \Leftrightarrow \! (d)$, p.~108) we obtain \eqref{nul_ev1} from  (\ref{dflkbv}). By \eqref{beurkhhh}, the other statement \eqref{nul_ev2} then follows from \eqref{nul_ev1} and \eqref{bd2}. \cqfd

\paragraph{Proof of Theorem \ref{main}: }
we define 
$$ H^{(p)} :=  \big( \frac{_{_1}}{^{^p}} C_{2 b_p s}  (  T_-( \rT_{p } , \bep_{p} ) )\big)_{s\in \bbR_+} \quad \textrm{and} \quad \widetilde{H}^{(p)} := \big( \frac{_{_1}}{^{^p}} C_{2 b_p s}  (  T_+( \rT_{p} , \bep_{ p}) )\big)_{s\in \bbR_+} \; .$$
From the definition of $(T_-( \rT_{p} , \bep_{ p}), T_+( \rT_{p} , \bep_{ p}))$, we note that 
\begin{equation}
\label{concaeffi}
\widetilde H^{(p)} \! \oplus \!  H^{(p)}=  \big( \frac{_{_1}}{^{^p}} C_{2 b_p s}  (   \rT_{ p} , \bep_{p}) \big)_{s\in \bbR_+} \; ;
\end{equation}
see also Figure \ref{fig: rooting}. 
Let $H$ and $\widetilde{H}$ be two independent  processes with the same law  $\bN (\,\cdot \, | \,\Gamma(H)\! = \! r)$. 
By Lemma \ref{contigu}, we deduce from Lemma \ref{zouloufli}, Lemma \ref{nupzer} and Proposition \ref{bdkjfvs} the following weak convergence on $\bC (\bbR_+, \bbR_+)\! \times \! \bC (\bbR_+, \bbR_+)\! \times \! \bbR_+\! \times \! \bbR_+$:
$$ \big(H^{(p)} ,   \widetilde{H}^{(p)} , \zeta (H^{(p)}), \zeta (\widetilde{H}^{(p)}) \big) \;  \overset{{\rm (d)}}{\underset{p\rightarrow \infty}{-\!\!\! -\!\!\! -\!\!\! \longrightarrow}} \; \big( H, \widetilde{H}, \zeta (H), \zeta( \widetilde{H}) \big) \; . $$
Then, along with (\ref{concaeffi}) and Lemma \ref{concatef}, we deduce from this the convergence \eqref{heigthcon} in  Theorem \ref{main}. \cqfd

\paragraph{Acknowledgement.} I am very grateful to Thomas Duquesne for discussions and suggestions, especially for his help on Proposition \ref{bdkjfvs}. Part of the work  was carried out during my visits at NYU Shanghai and at LaBRI. I thank both institutes for financial supports and my hosts Nicolas Broutin and Jean-Fran\c{c}ois Marckert for invitation.  I also acknowledge partial support from the Agence Nationale de la Recherche grant number ANR-14-CE25-0014 (ANR GRAAL).

\end{document}